\documentclass[preprint,12pt,3p]{elsarticle}
\usepackage{amssymb}
\usepackage{amsmath}
\usepackage{amsfonts}
\usepackage{comment}
\usepackage{here}
\usepackage{url}
\usepackage{color}
\usepackage{color}

\usepackage{amsthm}
\newtheorem{theorem}{Theorem}

\newcommand{\bl}[1]{\boldsymbol{#1}}
\newcommand{\mc}[1]{\mathcal{#1}}
\newcommand{\sign}{\text{sign}}
\usepackage {arydshln}

\journal{}
\begin{document}
	\begin{frontmatter}
		\title{
			Numerical verification for asymmetric solutions of the H\'enon equation on bounded domains 
		}
		\author[label1]{Taisei Asai \corref{cor1}}
		\address[label1]{Graduate School of Fundamental Science and Engineering, Waseda University, 3-4-1 Okubo, Shinjuku, Tokyo 169-8555, Japan}
		\cortext[cor1]{Corresponding author.}
		\ead{captino@fuji.waseda.jp}
		\author[label2]{Kazuaki Tanaka}
		\address[label2]{Institute for Mathematical Science, Waseda University, 3-4-1 Okubo, Shinjuku, Tokyo 169-8555, Japan}
		\ead{tanaka@ims.sci.waseda.ac.jp}
		\author[label5]{Shin'ichi Oishi}
		\address[label5]{Faculty of Science and Engineering, Waseda University, 3-4-1 Okubo, Shinjuku, Tokyo 169-8555, Japan}
		\ead{oishi@waseda.jp}
		\begin{abstract}
			The H\'enon equation, a generalized form of the Emden equation, admits symmetry-breaking bifurcation for a certain ratio of the transverse velocity to the radial velocity. Therefore, it has asymmetric solutions on a symmetric domain
			even though the Emden equation has no asymmetric unidirectional solution on such a domain.
			We discuss a numerical verification method for proving the existence of solutions of the H\'enon equation on a bounded domain.
			By applying the method to a line-segment domain and a square domain, we numerically prove the existence of solutions of the H\'enon equation for several parameters representing the ratio of transverse to radial velocity.
			As a result, we find a set of undiscovered solutions with three peaks on the square domain.
		\end{abstract}
		\begin{keyword}
			H\'enon equation \sep Numerical verification \sep Symmetry-breaking bifurcation \sep Elliptic boundary value problem
		\end{keyword}
	\end{frontmatter}
	\section{Introduction}
	The H\'enon equation was proposed as a model for mass distribution in spherically symmetric star clusters, which is important in studying the stability of rotating stars \cite{1}.
	One important aspect of the model is the Dirichlet boundary value problem
	\begin{align}
		\label{henon}
		\begin{cases}
			-\Delta u=|\boldsymbol{x}-\boldsymbol{x_0}|^{l}|u|^{p-1} u  & {\rm in} \quad \Omega,\\
			u=0  & {\rm on} \quad \partial \Omega,
		\end{cases}
	\end{align}
	where $\Omega \subset \mathbb{R}^{N}(N=1,2,3)$ is a bounded domain, $\boldsymbol{x}$ is the location of the star, and $u$ stands for the stellar density.
	Particularly, $\boldsymbol{x_0}$ is set at the center of the domain.
	The parameter $2 \leq p < p^*$ ($p^*=\infty$ if $N=1,2$ and $p^*=5$ if $N = 3$) is the polytropic index, determined according to the central density of each stellar type.
	The parameter $l \geq 0$ is the ratio of the transverse velocity to the radial velocity.
	These velocities can be derived by decomposing the space velocity vector into the radial and transverse components.
	
	When $l=0$, the H\'enon equation coincides with the Emden equation $-\Delta u=|u|^{p-1} u \text{  in  } \Omega$.
	In this case, the transverse velocity vanishes and the orbit becomes purely radial.
	Gidas, Ni, and Nirenberg proved that the Emden equation has no asymmetric unidirectional solution in a rectangle domain \cite{GNN}. However, Breuer, Plum, and McKenna reported some asymmetric solutions obtained with an approximate computation based on the Galerkin method \cite{EmdenI},
	which were called ``spurious approximate solutions'' caused by discretization errors.
	This example shows the need to verify approximate computations.
	By contrast, 
	a theoretical analysis \cite{SB} for large $l$ (when the orbit tends to be purely circular) found that the H\'enon equation admits symmetry-breaking bifurcation, thereby having several asymmetric solutions even on a symmetric domain.
	
	The importance of the H\'enon equation has led to active mathematical study on it over the last decade.
	For example, Amadori and Gladiali \cite{2} analyzed the bifurcation structure of \eqref{henon} with respect to parameter $p$.
	They applied an analytical method to the H\'enon equation that had worked for the Emden equation.
	Additionally, several numerical studies have been conducted on the H\'enon equation \cite{3,4,unitdisk,unitcub}.
	In particular, we are motivated by the work of Yang, Li, and Zhu \cite{3},
	who developed an effective computational method to find multiple asymmetric solutions of \eqref{henon} on the unit square $\Omega=(0,1)^2$ using algorithms based on the bifurcation method. They generated the bifurcation curve of \eqref{henon} with $p=3$ and numerically predicted bifurcation points around $l=0.5886933$ and $l=2.3654862$ using approximate computations.
	
	The purpose of our study is to prove the existence of solutions of \eqref{henon} using the Newton--Kantorovich theorem (see Theorem 2).
	We prove their existence through the following steps:
	\begin{enumerate}
		\item We construct approximate solutions $\hat{u}$ using the Galerkin method with polynomial approximations.
		\item Using the Newton--Kantorovich theorem (Theorem \ref{theo:nk}), we prove the existence of solutions $u$ of \eqref{henon} with nearby approximations $\hat{u}$ while sharply evaluating the error bound between $u$ and $\hat{u}$ in terms of the $H^1_0$-norm $\|\nabla \cdot\|_{L^2}$.
	\end{enumerate}
	By applying the above steps to the problem \eqref{henon} on the domains $\Omega=(0,1)^N$ $(N=1,2)$,
	we successfully prove the existence of several solutions for $l=0,2,4$.
	In particular, we find a set of solutions with three peaks, which were not revealed in \cite{3} (see Figure \ref{zu2}).
	
	The remainder of this paper is organized as follows.
	Some notation is introduced in Section \ref{sec:pre}.
	Sections \ref{sec:NV} and \ref{sec:eva} describe numerical verification based on the Newton--Kantorovich theorem together with evaluations of several required constants.
	Section \ref{sec:result} shows the results numerically proving the existence of several asymmetric solutions of \eqref{henon}.
	Subsequently, we discuss the solution curves of the problem for $p=3$ based on an approximate computation.
	\section{Preliminaries} \label{sec:pre}
	We begin by introducing some notation.
	For two Banach spaces $X$ and $Y$,
	the set of bounded linear operators from $X$ to $Y$ is denoted by $\mathcal{L} (X,Y)$.
	The norm of $T \in \mathcal{L} (X,Y)$ is defined by 
	\begin{align}
		\label{eq:normofdual}
		\|T\|_{\mathcal{L}(X, Y)}:=\sup _{0 \neq u \in X} \frac{\|T u\|_{Y}}{\|u\|_{X}} .
	\end{align}
	Let  $L^{p}(\Omega)$ $(1 \leq p<\infty)$ be the function space of $p$-th power Lebesgue integrable functions over a domain $\Omega$ with
	the $L^{p} $-norm
	$	\|u\|_{L^{p}}:=\left( \int_{\Omega} | u(x) |^p dx \right)^{1/p} < \infty.$
	When $p = 2$, $L^2 (\Omega)$ is the Hilbert space
	with the inner product
	$	(u, v)_{L^{2} }:=\int_{\Omega} u(x) v(x) d x$.
	Let $L^{\infty}(\Omega)$ be the function space of Lebesgue measurable functions over $\Omega$, with the norm $\|u\|_{L^{\infty}} :=\operatorname{ess} \sup \{|u(x)| : x \in \Omega\} \text{ for } u \in L^{\infty}(\Omega)$.
	We denote the first-order $L^2$ Sobolev space in $\Omega$ as $H^{1}(\Omega)$ and define
	\begin{align*}
		H_{0}^{1}(\Omega) :=\left\{u \in H^{1}(\Omega) : u=0 \text{ on } \partial \Omega \text{ in the trace sense}\right\}
	\end{align*}
	as the solution space for the target equation \eqref{henon}.
	We endow $H^1_0 (\Omega)$ with the inner product and norm
	\begin{align}
		(u,v)_{H^1_0}:&=(\nabla u,\nabla v)_{L^2}+\tau(u,v)_{L^2}, \quad u,v\in H^1_0(\Omega), \label{Vinner} \\
		\|u\|_{H^1_0}:&=\sqrt{(u,u)_{H^1_0}},\quad u\in H^1_0(\Omega), \label{Vnorm}
	\end{align}
	where $\tau$ is a nonnegative number chosen as \begin{align}
		\tau > -p|\bl{x}-\bl{x_0}|^{l}|\hat{u}(\bl{x})|^{p-1} \text{~~a.e.~~} \bl{x}\in \Omega \label{eq:tau}
	\end{align}
	for a numerically computed approximation $\hat{u} \in H^1_0(\Omega)$.
	The condition \eqref{eq:tau} is required in Subsection \ref{subsec:invnorm} and $\hat{u}$ is explicitly constructed in Section \ref{sec:result}.
	Because the norm $\| \cdot \|_{H^1_0}$ monotonically increases with respect to $\tau$, the $H^1_0 (\Omega)$ norm $\| \nabla \cdot \|_{L^2}$ is dominated by the norm $\| \cdot \|_{H^1_0}$ for all $\tau \geq 0$.
	Therefore, the error bound $ \|u-\hat{u}\|_{H^1_0} $ is always an upper bound for $\| \nabla (u-\hat{u}) \|_{L^2}$.
	The topological dual space of $H^1_0 (\Omega)$ is denoted by $H^{-1}$ with the usual supremum norm defined in \eqref{eq:normofdual}.
	
	The bound for the embedding $H^1_0(\Omega) \hookrightarrow L^{p}(\Omega)$ is denoted by $C_p$ $(p \geq 2)$. More precisely, $C_p$ is a positive number satisfying
	\begin{align}
		\|u\|_{L^{p}} \leq C_{p}\|u\|_{H^1_0} \quad \text{ for all } u \in H^1_0(\Omega). \label{umekomi}
	\end{align}	
	Note that $\|u\|_{H^{-1}} \leq C_{p}\|u\|_{L^{p'}}$, $u \in L^{p'}(\Omega)$ holds for $p'$ satisfying $p^{-1}+p'^{-1}=1$ .
	Explicitly estimating the embedding constant $C_p$ is important for our numerical verification.
		When $p=2$, we use the following optimal inequality:
	\begin{align*}
		\|u\|_{L^{2}} \leq \frac{1}{\sqrt{\lambda_{1}+\tau}}\|u\|_{H^1_0},
	\end{align*}
	where $\lambda_{1}$ is the first eigenvalue of the minus Laplacian in the weak sense.
	Especially when $\Omega=(0,1)^N$, we have $\lambda_{1}=N \pi^{2}$. 
	When $p$ is not $2$, we use the following theorems depending on the dimension of $\Omega$.
	We use \cite[Lemma 7.12]{nakao2019book} to obtain an explicit value of $C_p$ for a one-dimensional bounded domain.\\
	\begin{theorem}[{\cite[Lemma 7.12]{nakao2019book}}]
	    \label{theo:1D_Cp}
	    Let $\Omega=(a, b) \subset \mathbb{R}$, with $a \in \mathbb{R} \cup\{-\infty\}$, $ b \in \mathbb{R} \cup\{+\infty\}$, $a<b$.
	    Moreover, let $\rho ^*$ denote the minimal point of the spectrum of $-u''$ on $H^1_0 (\Omega)$, i.e. $\rho^{*}=\pi^{2} /(b-a)^{2}$ if $(a,b)$ is bounded.
	    Then, for all $u \in H_{0}^{1}(\Omega)$,
	    \begin{align*}
            \|u\|_{L^{p}} \leq C_{p}\|u\|_{H_{0}^{1}} \quad(p \in(2, \infty)),
        \end{align*}
        where, abbreviating $\varepsilon:=\frac{2}{p} \in(0,1)$,
        \begin{align*}
            C_{p}:=
            \begin{cases}
			\frac{1}{\sqrt{2}}(1-\varepsilon)^{\frac{1}{4}(1-\varepsilon)}(1+\varepsilon)^{\frac{1}{4}(1+\varepsilon)} \tau^{-\frac{1}{4}(1+\varepsilon)} &\text { if } \quad \rho^{*} \leq \tau \frac{1-\varepsilon}{1+\varepsilon},\\
			\frac{1}{\sqrt{\rho^{*}+\tau}}\left(\rho^{*}\right)^{\frac{1}{4}(1-\varepsilon)} &\text { otherwise },
		\end{cases}
        \end{align*}
        for $p \in(2, \infty)$.
	\end{theorem}

    When $N \geq 2$, we use \cite[Corollary A.2]{tanakaumekomi} or \cite[Lemma 7.10]{nakao2019book} to obtain $C_p$ for bounded domains $\Omega \subset \mathbb{R}^{N}$.
	In our numerical experiments in Section \ref{sec:result}, $C_p$ evaluated by \cite[Corollary A.2]{tanakaumekomi} is smaller than that evaluated by \cite[Lemma 7.10]{nakao2019book}.
	
	\begin{theorem}[{\cite[Corollary A.2]{tanakaumekomi}}]
		\label{sitate}
		Let $\Omega \subset \mathbb{R}^{N}(N \geq 2)$ be a bounded domain, the measure of which is denoted by $| \Omega |$.
		Let $p \in (N/(N-1),2N/(N-2) ]$ if $N \geq 3$, $p \in (2, \infty)$ if $N=2$.
		We set $q=Np/(N+p)$. Then, \eqref{umekomi} holds for
		\begin{align*}
			C_p( \Omega ) = | \Omega |^{\frac{2-q}{2q}}T_p.
		\end{align*}
		Here, $T_p$ is defined by
		\begin{align*}
			T_p=\pi^{-\frac{1}{2}} N^{-\frac{1}{q}}\left(\frac{q-1}{N-q}\right)^{1-\frac{1}{q}}\left\{\frac{\Gamma\left(1+\frac{N}{2}\right) \Gamma(N)}{\Gamma\left(\frac{N}{q}\right) \Gamma\left(1+N-\frac{N}{q}\right)}\right\}^{\frac{1}{N}},
		\end{align*}
		where $\Gamma$ is the gamma function.
	\end{theorem}

	\section{Numerical verification method} \label{sec:NV}
	This section discusses the numerical verification method used in this paper.
	We first define the operator $f$ as
	\begin{align*}
		f:
		\begin{cases}
			u(\cdot ) &\mapsto  |\boldsymbol{\cdot}-\boldsymbol{x_0}|^{l}|u(\cdot)|^{p-1} u(\cdot),\\
			H^1_0 (\Omega )&\to H^{-1},
		\end{cases}
	\end{align*}
	where $2 \leq p < p^*$ ($p^*=\infty$ if $N=1,2$ and $p^*=5$ if $N = 3$).
	Furthermore, we define the nonlinear operator $F:H^1_0(\Omega) \rightarrow H^{-1}$ by $F(u):=-\Delta u -f(u)$, which is given by
	\begin{align}
		\langle F(u),v\rangle = (\nabla u,\nabla v)_{L^2} -\langle f(u),v \rangle \quad \text{for all} \quad v\in H^1_0(\Omega), \nonumber
	\end{align}
	where $\langle f(u), v\rangle=\int_{\Omega} (|\bl{x}-\bl{x_0}|^{l} |u(\bl{x})|^{p-1} u(\bl{x})) v(\bl{x}) d \bl{x}$.
	The Fr\'echet derivatives of $f$ and $F$ at $\varphi \in H^1_0(\Omega)$ are denoted by $f'_{\varphi}$ and $F'_{\varphi}$, respectively, and given by
	\begin{align}
		\left\langle f_{\varphi}^{\prime} u, v\right\rangle &=\int_{\Omega} (p |\bl{x}-\bl{x_0}|^{l} |\varphi(\bl{x})|^{p-1}) u(\bl{x}) v(\bl{x}) d \bl{x} \quad \text{ for all } \quad u, v \in H^1_0(\Omega) ,\label{eq:fd_dual}\\
		\left\langle {F}_{\varphi}^{\prime} u, v\right\rangle &=(\nabla u, \nabla v)_{L^{2}}-\left\langle f_{\varphi}^{\prime} u, v\right\rangle \quad \text{ for all } \quad u, v \in H^1_0(\Omega) .
	\end{align}
	Then, we consider the following problem:
	\begin{align}
		\label{weakform}
		\text{Find} \quad u \in H^1_0(\Omega) \quad \text{s.t.} \quad F(u)=0,
	\end{align}
	which is the weak form of the problem \eqref{henon}.
	To conduct the numerical verification for this problem, we apply the Newton--Kantorovich theorem,
	which enables us to prove the existence of a true solution $u$ near a numerically computed ``good'' approximate solution $\hat{u}$ (see, for example, \cite{5}).
	Hereafter, $B(\hat{u},r)$ and $\bar{B}(\hat{u},r)$ respectively denote the open and closed balls with center approximate solution $\hat{u}$ and radius $r$ in terms of norm $\| \cdot \|_{H^1_0}$.
	\begin{theorem}[Newton--Kantorovich's theorem]
		\label{theo:nk}
		Let $\hat{u} \in H^1_0 (\Omega)$ be some approximate solution of $F(u)=0$. 
		Suppose that there exists some $\alpha >0$ satisfying
		\begin{align}
			\label{NKalpha}
			\| F'^{-1}_{\hat{u}}F({\hat{u}})\| _{H^1_0} \leq \alpha.
		\end{align}
		Moreover, suppose that there exists some $\beta >0$ satisfying
		\begin{align}
			\label{NKbeta}
			\|F'^{-1}_{\hat{u}}(F'_v-F'_w)\| _{\mathcal{L}({H^1_0},{H^1_0})} \leq \beta \|v-w\|_{H^1_0} , \quad \text{ for all } v,w \in D, 
		\end{align}
		where
		$D=B(\hat{u},2 \alpha + \delta)$ is an open ball depending on the above value $\alpha >0$ for small $\delta>0$.
		If 
		\begin{align*}
			\alpha \beta \leq \frac{1}{2},
		\end{align*}
		then there exists a solution $u \in H^1_0(\Omega)$ of $F(u) = 0 $ in $\bar{B}(\hat{u},\rho)$ with
		\begin{align}
			\rho =\frac{1-\sqrt{1-2\alpha \beta}}{\beta}.\nonumber
		\end{align}
		Furthermore, 
		the solution $u$ is unique in $\bar{B}(\hat{u},2\alpha)$.
	\end{theorem}
	\section{Evaluation for $\alpha$ and $\beta$} \label{sec:eva}
	To apply Theorem \ref{theo:nk} to the numerical verification for problem \eqref{henon},
	we need to explicitly evaluate $\alpha$ and $\beta$.
	The left side of \eqref{NKalpha} is evaluated as
	\begin{align*}
		\left\|F_{\hat{u}}^{\prime-1} F(\hat{u})\right\|_{H^1_0} \leq\left\|F'^{-1}_{\hat{u}}\right\|_{\mathcal{L}\left(H^{-1}, {H^1_0}\right)}\|F(\hat{u})\|_{H^{-1}}.
	\end{align*}
	Therefore, we set
	\begin{align*}
		\alpha = \left\|F'^{-1}_{\hat{u}}\right\|_{\mathcal{L}\left(H^{-1}, {H^1_0}\right)}\|F(\hat{u})\|_{H^{-1}}.
	\end{align*}
	Moreover, the left side of (\ref{NKbeta}) is estimated as
	\begin{align*}
		\left\|F_{\hat{u}}^{\prime-1}\left(F_{v}^{\prime}-F_{w}^{\prime}\right)\right\|_{\mathcal{L}({H^1_0}, {H^1_0})} &\leq\left\|F_{\hat{u}}^{\prime-1}\right\|_{\mathcal{L}\left(H^{-1}, {H^1_0}\right)}\left\|F_{v}^{\prime}-F_{w}^{\prime}\right\|_{\mathcal{L}\left({H^1_0}, H^{-1}\right)}\\
		&=\left\|F_{\hat{u}}^{\prime-1}\right\|_{\mathcal{L}\left(H^{-1}, {H^1_0}\right)}\left\|f'_v-f'_w\right\|_{\mathcal{L}\left({H^1_0}, H^{-1}\right)}.
	\end{align*}
	Hence, the desired value of $\beta$ is obtained via 
	\begin{align*}
		\beta &\leq \| F'^{-1}_{\hat{u}}\| _{\mathcal{L}({H^{-1}},{H^1_0})} L,
	\end{align*}
	where $L$ is the Lipschitz constant satisfying 
	\begin{align}
		\left\|f'_v-f'_w\right\|_{\mathcal{L}\left({H^1_0}, H^{-1}\right)}
		\leq L\|v-w\|_{H^1_0} \quad \text{ for all } v, w \in D . \label{Lip-satis}
	\end{align} 
	We are left to evaluate the inverse operator norm $\| F'^{-1}_{\hat{u}}\|_{\mathcal{L}({H^{-1}},{H^1_0})}$, the residual norm $\|F({\hat{u}})\| _{H^{-1}}$, and the Lipschitz constant $L$ for problem \eqref{weakform}.
	\subsection{Residual norm $\|F({\hat{u}})\| _{H^{-1}}$}
	If the approximation $\hat{u}$ is sufficiently smooth so that $\Delta \hat{u} \in L^2 (\Omega)$, we can evaluate the residual norm $\|F({\hat{u}})\| _{H^{-1}}$ as follows:
	\begin{align}
		\label{eq:res}
		\|F({\hat{u}})\| _{H^{-1}} \leq C_{2}\|\Delta \hat{u}+f(\hat{u})\|_{L^{2}},
	\end{align}
	where $C_2$ is the embedding constant satisfying \eqref{umekomi} for $p=p'=2$.
	Our numerical experiments discussed in Section \ref{sec:result} use this evaluation, calculating the $L^2$-norm via stable numerical integration with all rounding errors strictly estimated.
	
	However, the condition $\Delta \hat{u} \in L^2 (\Omega)$ is not satisfied such as when we construct $\hat{u}$ with a piecewise linear finite element basis.
	We use the method of \cite[Subsection 7.2]{nakao2019book} to evaluate the residual norm applicable to such a case.
	The following is a brief description of the evaluation method.
	First, we find an approximation $\rho \in H(\operatorname{div}, \Omega)=\left\{\tau \in L^{2}(\Omega)^{N} : \operatorname{div} \tau \in L^{2}(\Omega)\right\}$ to $\nabla \hat{u}$.
	Then, the residual norm is evaluated as 
	\begin{align*}
		\|{F}(\hat{u})\|_{H^{-1}} &=\|-\Delta \hat{u}-f(\hat{u})\|_{H^{-1}}, \\
		&=\|-\Delta \hat{u}+\operatorname{div} \rho-\operatorname{div} \rho-f(\hat{u})\|_{H^{-1}}, \\
		& \leq\|\operatorname{div}(-\nabla \hat{u}+\rho)\|_{H^{-1}}+\|\operatorname{div} \rho+f(\hat{u})\|_{H^{-1}}, \\
		& \leq\|-\nabla \hat{u}+\rho\|_{L^{2}}+C_{2}\|\operatorname{div} \rho+f(\hat{u})\|_{L^{2}},
	\end{align*}
	where we used $\|\mathrm{div} \omega\|_{H^{-1}} \leq\|\omega\|_{L^{2}}$ for $\omega \in H(\mathrm{div}, \Omega)$.
	As mentioned in \cite[Subsection 7.2]{nakao2019book}, $\rho$ can be computed without additional computational resources when we use the mixed finite element method to construct $\hat{u}$.
	\subsection{Inverse operator norm $\| F'^{-1}_{\hat{u}}\|_{\mathcal{L}({H^{-1}},{H^1_0})}$} \label{subsec:invnorm}
	In this subsection, we evaluate the inverse operator norm $\| F'^{-1}_{\hat{u}}\|_{\mc{L}(H^{-1},H^1_0)}$.
	To this end, we use the following theorem.
	\begin{theorem}[\cite{plum2009computer}]\label{invtheo}
		Let $\Phi:{H^1_0}(\Omega )\rightarrow H^{-1}$ be the canonical isometric isomorphism; that is, $\Phi$ is given by
		\begin{align*}
			\left\langle\Phi u,v\right\rangle:=\left(u,v\right)_{H^1_0}
			~~~{\rm for}~u,v\in {H^1_0(\Omega)}.
		\end{align*}
		If
		\begin{align}
			\displaystyle \mu_{0}
			:=\min\left\{|\mu|\ :\ \mu\in\sigma_{p}\left(\Phi^{-1}{F}_{\hat{u}}'\right)\cup\{1\}\right\}
			\label{mu0}
		\end{align}
		is positive,
		then the inverse of ${F}_{\hat{u}}'$ exists, and we have
		\begin{align}
			\left\|{F}_{\hat{u}}^{\prime-1}\right\|_{\mathcal{L}(H^{-1},H^1_0)}\leq\mu_{0}^{-1},\label{Ktheo}
		\end{align}
		where $\sigma_{p}\left(\Phi^{-1}{F}_{\hat{u}}'\right)$ denotes the point spectrum of $\Phi^{-1}{F}_{\hat{u}}'$.
	\end{theorem}

	The eigenvalue problem $\Phi^{-1}{F}_{\hat{u}}'u=\mu u$ in ${H^1_0(\Omega )}$ is equivalent to
	\begin{align}
		\label{eq:eigforTable}
		\left(\nabla u,\nabla v\right)_{L^2}-\langle f'_{\hat{u}} u,v\rangle
		=\mu\left( u,v\right)_{H^1_0}~~{\rm for~all}~v\in {H^1_0(\Omega )},
	\end{align}
	where $ \left( u,v\right)_{H^1_0} $ denotes the inner product defined in \eqref{Vinner} that depends on $\tau$ and $\langle f'_{\hat{u}} u,v\rangle$ is given by \eqref{eq:fd_dual}.

    We consider the operator $\mathcal{N}:=\Phi -{F}_{\hat{u}}^{\prime}$ from $H^1_0 (\Omega)$ to $H^{-1}$, which satisfies $\langle \mathcal{N} u, v \rangle =\int_{\Omega} (p |\bl{x}-\bl{x_0}|^{l} |\hat{u}(\bl{x})|^{p-1}) u(\bl{x}) v(\bl{x}) d \bl{x}$ for all $u,v \in H^1_0 (\Omega)$.
	Because $\mathcal{N}$ maps $H^1_0 (\Omega)$ into $L^2 (\Omega)$ and the embedding $L^2 (\Omega) \hookrightarrow H^{-1}$ is compact,
	$\mathcal{N}:H^1_0 (\Omega) \to H^{-1}$ is a compact operator.
	Therefore, $F'_{\hat{u}}$ is a Fredholm operator, and the spectrum $\sigma \left(\Phi^{-1}{F}_{\hat{u}}'\right)$ of $\Phi^{-1}{F}_{\hat{u}}'$ is given by 
	\begin{align*}
    \sigma\left(\Phi^{-1} F_{\hat{u}}^{\prime}\right)=1-\sigma\left(\Phi^{-1} \mathcal{N}\right)=1-\left\{\sigma_{p}\left(\Phi^{-1} \mathcal{N}\right) \cup\{0\}\right\}=\sigma_{p}\left(\Phi^{-1} F_{\hat{u}}^{\prime}\right) \cup\{1\}.
    \end{align*}
	Accordingly, it suffices to look for eigenvalues $\mu\neq 1$.
	By setting $\lambda=(1-\mu)^{-1}$, we further transform this eigenvalue problem into
	\begin{align}
		{\rm Find}~u\in H^1_0(\Omega)~{\rm and}~\lambda\in \mathbb{R}~{\rm s.t.}~
		(u,v)_{H^1_0}
		=\lambda\langle(\tau+f'_{\hat{u}})u,v\rangle~~{\rm for~all}~v\in {H^1_0(\Omega )},\label{eiglam}
	\end{align}
	where $\langle(\tau+f'_{\hat{u}})u,v\rangle= \int_{\Omega} (\tau + p |\bl{x}-\bl{x_0}|^{l} |\hat{u}(\bl{x})|^{p-1}) u(\bl{x}) v(\bl{x}) d \bl{x}$ for $u,v \in H^1_0 (\Omega)$.
	Because $\tau$ is chosen so that $\tau+f'_{\hat{u}}$ becomes positive (see \eqref{eq:tau}), \eqref{eiglam} is a regular eigenvalue problem, the spectrum of which consists of a sequence $\{\lambda_{k}\}_{k=1}^{\infty}$ of eigenvalues converging to $+\infty$.
	To compute $\| F'^{-1}_{\hat{u}}\|_{\mathcal{L}({H^{-1}},{H^1_0})}$ on the basis of Theorem \ref{invtheo},
	we need to enclose the eigenvalue $\lambda$ of \eqref{eiglam} that minimizes the corresponding absolute value of $|\mu|\left(=|1-\lambda^{-1}|\right)$.
	We consider the approximate eigenvalue problem
	\begin{align}
		{\rm Find}~u_{M} \in V_{M}~{\rm and}~\lambda^{M}\in \mathbb{R}~{\rm s.t.}
		\left(u_{M}, v_{M}\right)_{H^1_0}
		=\lambda^{M}\langle(\tau+f'_{\hat{u}})u_{M},v_{M}\rangle~~{\rm for~all}~v_{M}\in V_{M},\label{applam}
	\end{align}
	where $V_{M}$ is a finite-dimensional subspace of ${H^1_0(\Omega )}$ such as the space spanned by the finite element basis and Fourier basis.
	For our problem, $V_{M}$ will be explicitly chosen in Section \ref{sec:result}.
	Note that \eqref{applam} is a matrix problem
	with eigenvalues that can be enclosed with rigorous computation techniques (see, for example, \cite{behnke1991calculation,rump1999book,miyajima2012numerical}).
	
	We then estimate the error between the $k$-th eigenvalue $\lambda_{k}$ of \eqref{eiglam} and the $k$-th eigenvalue $\lambda_{k}^{M}$ of \eqref{applam}.
	We consider the weak formulation of the Poisson equation,
	\begin{align}
		\left(u_g, v\right)_{H^1_0}=\left(g,v\right)_{L^2}~~~{\rm for~all}~v\in H^1_0 (\Omega)\label{poisson}
	\end{align}
	given $g\in L^{2}\left(\Omega\right)$.
	This equation has a unique solution $u_g\in {H^1_0(\Omega )}$ for each $g\in L^{2}\left(\Omega\right)$ \cite{grisvard2011elliptic}.
	Let $P^{\tau}_{M}:{H^1_0(\Omega )} \rightarrow V_{M}$ be the orthogonal projection defined by
	\begin{align*}
		\left(P^{\tau}_{M}u-u,v_{M}\right)_{H^1_0}=0~~~{\rm for~all}~u\in H^1_0(\Omega){\rm~and~}v_{M}\in V_{M}.
	\end{align*}
	The following theorem enables us to estimate the error between $\lambda_{k}$ and $\lambda_{k}^{M}$.
	\begin{theorem}[\cite{tanaka2014verified, liu2015framework}]\label{eigtheo}
		Let $\hat{u} \in H^1_0(\Omega) \cap L^{\infty}(\Omega)$.
		Suppose that there exists $C^{\tau}_{M}>0$ such that
		\begin{align}
			\left\|u_{g}-P^{\tau}_{M}u_{g}\right\|_{H^1_0}\leq C^{\tau}_{M}\left\|g\right\|_{L^{2}}\label{CN}
		\end{align}
		for any $g\in L^{2}\left(\Omega\right)$ and the corresponding solution $u_{g}\in {H^1_0(\Omega )}$ of \eqref{poisson}.
		Then,
		\begin{align*}
			\displaystyle \frac{\lambda_{k}^{M}}{\lambda_{k}^{M}(C^{\tau}_{M})^{2}\|\tau+f'_{\hat{u}} \|_{L^{\infty}}+1}\leq\lambda_{k}\leq\lambda_{k}^{M},
		\end{align*}
		where the $L^{\infty}$-norm is defined by $\|\tau+f'_{\hat{u}} \|_{L^{\infty}}:= \text{\rm esssup\,} \{|\tau+p|\bl{x}-\bl{x_0}|^{l}|\hat{u}(\bl{x})|^{p-1}|\,:\,\bl{x} \in \Omega\}$.
	\end{theorem}
	The right inequality is known as the Rayleigh--Ritz bound, which is derived from the min-max principle:
	\begin{align*}
		\displaystyle \lambda_{k}=\min_{H_{k}\subset {H^1_0(\Omega )}}\left(\max_{v\in H_{k}\backslash\{0\}}\frac{\left\| v\right\|_{H^1_0}^{2}}{\left\|av\right\|_{L^{2}}^{2}}\right)\leq\lambda_{k}^{M},
	\end{align*}
	where $a(\bl{x})=\sqrt{\tau+p|\boldsymbol{x}-\boldsymbol{x_0}|^{l}|\hat{u}(\bl{x})|^{p-1}}$,
	and the minimum is taken over all $k$-dimensional subspaces $H_{k}$ of ${H^1_0(\Omega )}$.
	The left inequality was proved in \cite{tanaka2014verified, liu2015framework}.
	Assuming the $H^{2}$-regularity of solutions to \eqref{poisson} (which follows, for example, when $\Omega$ is a convex polygonal domain {\rm \cite[Section 3.3]{grisvard2011elliptic}}), {\rm \cite[Theorem 4]{tanaka2014verified}} ensures the left inequality.
	A more general statement that does not require the $H^{2}$-regularity is proved in \cite[Theorem 2.1]{liu2015framework}.
	
	When the solution of \eqref{poisson} has $H^{2}$-regularity, \eqref{CN} can be replaced with
	\begin{align}
		\label{eq:CNtau}
		\left\|u-P^{\tau}_{M}u\right\|_{H^1_0}\leq C^{\tau}_{M}\left\|-\Delta u + \tau u\right\|_{L^{2}}~~~{\rm for~all}~u\in H^{2}(\Omega)\cap H^1_0(\Omega).
	\end{align}
	The constant $C^{\tau}_{M}$ satisfying \eqref{eq:CNtau} is obtained as
	$C_{M}^{\tau}=C_{M} \sqrt{1+\tau\left(C_{M}\right)^{2}}$ (see \cite[Remark A.4]{tanaka2017numerical}), where we denote $C_{M}=C_{M}^{0}$ with $\tau=0$.
	For example, when $\Omega = (0, 1)^N$, an explicit value of $C_{M}$ is obtained for $V_{M}$ spanned by the Legendre polynomial basis using {\rm \cite[Theorem 2.3]{6}}.
	This will be used for our computation in Section \ref{sec:result}.
	\begin{theorem}[\cite{6}]\label{hokangosatheo}
		When $\Omega = (0,1)^N$, the inequality
		\begin{align*}
			\left\|\nabla\left(u-P_{M} u\right)\right\|_{L^{2}} \leq C_{M}\|\Delta u\|_{L^{2}}
			~~
			\text{for~all}
			~
			u \in H^{2}(\Omega)\cap H^1_0(\Omega)
		\end{align*}
		holds for 
		\begin{align*}
			&C_{M}=\max \left\{\frac{1}{2(2 M+1)(2 M+5)}+\frac{1}{4(2 M+5) \sqrt{2 M+3} \sqrt{2 M+7}}\right. , \\
			&\left.\frac{1}{4(2 M+5) \sqrt{2 M+3} \sqrt{2 M+7}}+\frac{1}{2(2 M+5)(2 M+9)}+\frac{1}{4(2 M+9) \sqrt{2 M+7} \sqrt{2 M+11}}\right\}^{\frac{1}{2}}.
		\end{align*}
	\end{theorem}
	\subsection{Lipschitz constant $L$}
	Hereafter, we denote $d\,(=d(\Omega,l)):=\max\{ | \bl{x}-\bl{x_0} | ^l :\bl{x}  \in \Omega\}$.
	The Lipschitz constant $ L $ satisfying \eqref{Lip-satis}, which is required for obtaining $ \beta $, is estimated as follows:
	\begin{align}
		\label{eq:lip1}
		\left \|f'_v-f'_w\right\|_{\mathcal{L}\left({H^1_0}, H^{-1}\right)} \nonumber 
		&\leq p \sup_{0\neq \phi \in {H^1_0}}\sup_{0\neq \psi \in {H^1_0}}\frac { | \int_{\Omega}|\boldsymbol{x}-\boldsymbol{x_0}|^{l}(|v(\bl{x})|^{p-1}\phi(\bl{x}) -|w(\bl{x})|^{p-1}\phi(\bl{x}))  \psi(\bl{x}) d \bl{x} |}{\| \phi \|_{H^1_0} \| \psi \|_{H^1_0}}\nonumber \\
		&\leq p d \sup_{0\neq \phi \in {H^1_0}}\sup_{0\neq \psi \in {H^1_0}}\frac { | \int_{\Omega}(|v(\bl{x})|^{p-1} -|w(\bl{x})|^{p-1} ) \phi(\bl{x})  \psi(\bl{x}) d \bl{x}|}{\| \phi \|_{H^1_0} \| \psi \|_{H^1_0}}.
	\end{align}
 	Using the mean value theorem,
	the numerator of \eqref{eq:lip1} is evaluated as
	\begin{align*}
		&\left| \int_{\Omega}(|v(\bl{x})|^{p-1}-|w(\bl{x})|^{p-1})\phi(\bl{x}) \psi(\bl{x}) d \bl{x} \right| \\
		&=\biggl| \int_{\Omega}\int_{0}^{1}(p-1)\sign(w(\bl{x})+t(v(\bl{x})-w(\bl{x})))|w(\bl{x})+t(v(\bl{x})-w(\bl{x}))|^{p-2}  dt \biggl.\\
    		&~~~~~~~~~~~~~~~~~~~~~~~~~~~~~~~~~~~~~~~~~~~~~~~~~~~~~~~~~~~~~~~~~~\biggl.(v(\bl{x})-w(\bl{x}))\phi(\bl{x}) \psi(\bl{x}) d \bl{x} \biggl|\\
		&= (p-1)\biggl| \int_{0}^{1} \int_{\Omega} \sign(w(\bl{x})+t(v(\bl{x})-w(\bl{x})))|w(\bl{x})+t(v(\bl{x})-w(\bl{x}))|^{p-2}   \biggl.\\
		&~~~~~~~~~~~~~~~~~~~~~~~~~~~~~~~~~~~~~~~~~~~~~~~~~~~~~~~~~~~~~~~~~~\biggl. (v(\bl{x})-w(\bl{x})) \phi(\bl{x}) \psi(\bl{x}) d \bl{x} dt \biggl| \\
		& \leq (p-1)\int^{1}_{0} \|tv+(1-t)w\|^{p-2}_{L^{p+1}} \| v-w \|_{L^{p+1}}\| \phi \| _{L^{p+1} }\| \psi \| _{L^{p+1} } dt\\
		& \leq  (p-1)C^3_{p+1}\int^{1}_{0} \|tv+(1-t)w\|^{p-2}_{L^{p+1}} dt \| v-w \|_{H^1_0}\| \phi \| _{H^1_0}\| \psi \| _{H^1_0} \\
		& \leq  (p-1)C^3_{p+1}\max \left\{ \| v \|_{L^{p+1}}, \|w\|_{L^{p+1}}\right\}^{p-2} \| v-w \|_{H^1_0}\| \phi \| _{H^1_0}\| \psi \| _{H^1_0},
	\end{align*}
	for all $0\neq \phi , \psi \in H^1_0 (\Omega)$.
	Therefore, we have
	\begin{align}
		L \leq p (p-1) d C^3_{p+1} \max \left\{ \| v \|_{L^{p+1}}, \|w\|_{L^{p+1}}\right\}^{p-2}. \nonumber
	\end{align}
	Choosing  $v,w$ from $D=B(\hat{u},r)$, $r=2\alpha + \delta$ for small $\delta>0$, we can express them as
	\begin{align}
		\begin{cases}
			v=\hat{u}+r\eta ,  & \| \eta \|_{H^1_0} \leq 1,\\
			w=\hat{u}+r\xi , & \| \xi \|_{H^1_0} \leq 1 .\nonumber
		\end{cases}
	\end{align}
	Hence, it follows that
	\begin{align}
			L & \leq p (p-1) d C^3_{p+1} \max \left\{ \| \hat{u}+r\eta \|_{L^{p+1}}, \|\hat{u}+r\xi\|_{L^{p+1}}\right\}^{p-2} \nonumber \\
		    & \leq p (p-1) d C^3_{p+1} (\| \hat{u} \| _{L^{p+1} }+C_{p+1}r)^{p-2}. \label{eq:Lip_end}
	\end{align}

	\section{Numerical results}\label{sec:result}
	In this section, 
	we present numerical results where the existence of solutions of \eqref{henon} was proved for $p=3$ on the domains $\Omega=(0,1)^N ~(N=1,2) $ via the method presented in Sections \ref{sec:NV} and \ref{sec:eva}.
	All computations were implemented on a computer with 2.20 GHz Intel Xeon E7-4830 CPUs $\times$ 4, 2 TB RAM, and CentOS 7 using MATLAB 2019b with GCC Version 6.3.0.
	All rounding errors were strictly estimated using the toolboxes kv Library \cite{kv} Version 0.4.49 and Intlab Version 11 \cite{rump1999book}.
	Therefore, the accuracy of all results was guaranteed mathematically. 
	We constructed approximate solutions of \eqref{henon} from a Legendre polynomial basis discussed in \cite{6}.
	Specifically, we constructed approximate solutions $\hat{u}$ using the basis functions $\phi_{n}~(n=1,2,3, \cdots)$ defined by
	\begin{align}
		&\phi_{n}(x)=\frac{1}{n(n+1)} x(1-x) \frac{d Q_{n}}{d x}(x) \nonumber\\
		&~~~~~~~~~~~~~~~~\text{ with } Q_{n} (x)=\frac{(-1)^{n}}{n !}\left(\frac{d}{d x}\right)^{n} x^{n}(1-x)^{n}, \quad n=1,2,3, \cdots.
		\label{eq:Legendrepoly}
	\end{align}
	
	\subsection{Numerical results on the unit line-segment}\label{sec:1Dresult}
	To apply our method to $\Omega=(0,1)$,
    we define the finite-dimensional subspace $V_{M}$ of $H^1_0(\Omega)$ as
    \begin{align*}
		V_{M} &:= \left\{\sum_{i=1}^{M}  ~u_{i} \phi_{i}(x) ~:~ u_{i} \in \mathbb{R} \right\},
	\end{align*}
	where $2 \leq M< \infty$.
	We computed approximate solutions $\hat{u} \in V_{M}$ by solving the problem of the matrix equation
	\begin{align}
		{\rm Find}&~\hat{u} \in V_{M}~{\rm s.t.}~	
		\left(\nabla \hat{u},\nabla v_{M}\right)_{L^2}
		=\left(f(\hat{u}),v_{M}\right)_{L^2}~~{\rm for~all}~v_{M}\in V_{M} \label{eq:matrixeq}
	\end{align}
	using the usual Newton method.
	When we look for a symmetric solution,
	we restrict the solution space and its finite-dimensional subspace.
	The following subspace $V^1$ of $H^1_0(\Omega)$ is endowed with the same topology
	\begin{align}
		V^1 &:= \left\{u \in H^1_0(\Omega): u \text{ is symmetric with respect to } x=\frac{1}{2}\right\}. \label{eq:V1_def_1D}
	\end{align}
	Then, we define the finite-dimensional subspace $V^{1}_{M}~(M\geq 2)$ of $V^{1}$ as 
	\begin{align*}
		V^1_{M} &:= \left\{\sum_{i=1\atop \text{$i$ is odd}}^{M}  ~u_{i} \phi_{i}(x) ~:~ u_{i} \in \mathbb{R} \right\}.
	\end{align*}
	The method presented in Sections \ref{sec:NV} and \ref{sec:eva} can be directly applied when the function spaces $H^{1}_{0}(\Omega )$ and $V_{M}$ are replaced with $V^1$ and $V^{1}_M$, respectively.
	This restriction reduces the amount of calculation because 
	the matrices in \eqref{eq:matrixeq} become smaller.
	Moreover, because eigenfunctions of \eqref{applam} are also restricted to be symmetric, eigenvalues associated with anti-symmetric eigenfunctions drop out of the minimization in \eqref{mu0}.
	Therefore, the constant $K$ can be reduced. 
	The other constants required in the verification process (that is, $C_p$ and $\|F({\hat{u}})\| _{H^{-1}}$) are not affected by the restriction.
	Using the evaluation \eqref{eq:Lip_end} when $p = 3$ and $\Omega=(0,1)$ with the center $\boldsymbol{x_0}=(1/2)$, we evaluated the Lipschitz constant $L$ as
	\begin{align*}
		L \leq 6 \left( \frac{1}{2} \right) ^l   C^3_{4}(\| \hat{u} \| _{L^{4}}+C_{4}r).
	\end{align*}

	Table \ref{1D} shows the approximate solutions together with their verification results on $\Omega=(0,1)$.
	The red dashed lines indicate the symmetry of each solution.
	To satisfy inequality \eqref{eq:tau},
	our program set $\tau$ to the next floating-point number after a computed upper bound of the right side of \eqref{eq:tau}.
	Therefore, when $\hat{u}$ vanishes at some point on $\overline{\Omega}$,
	$\tau$ is set to the floating-point number after zero, which is approximately $4.9407\times 10^{-324}$.
	In Table \ref{1D}, $\|F(\hat{u})\|_{H^{-1}}$, $\|F'^{-1}_{\hat{u}}\|_{\mathcal{L}(H^{-1},H^1_0)}$, $L$, $\alpha$, and $\beta$ denote the constants required by Theorem \ref{theo:nk}.
	Moreover, $r_A$ and $r_R$ denote an upper bound for absolute error $\| u- \hat{u} \|_{H^1_0}$ and relative error $\| u- \hat{u} \|_{H^1_0}/\| \hat{u} \|_{H^1_0}$, respectively.
	The values in row ``Peak'' represent upper bounds for the maximum values of the corresponding approximate solutions in decimal form.
	
	The values in rows $\mu_{1}$--$\mu_{5}$ represent approximations of the five smallest eigenvalues of \eqref{eq:eigforTable} discretized in $V_{40} \subset H^1_0(\Omega)$, which is spanned by the basis functions $\phi_n$ ($n=1,2,\cdots,40$) without the  restriction of symmetry.
	When $l=2,4$, symmetric solutions have two negative eigenvalues and asymmetric solutions have one negative eigenvalue.
	
\begin{figure}[H]
		\centering
		\includegraphics[width=0.8\textwidth]{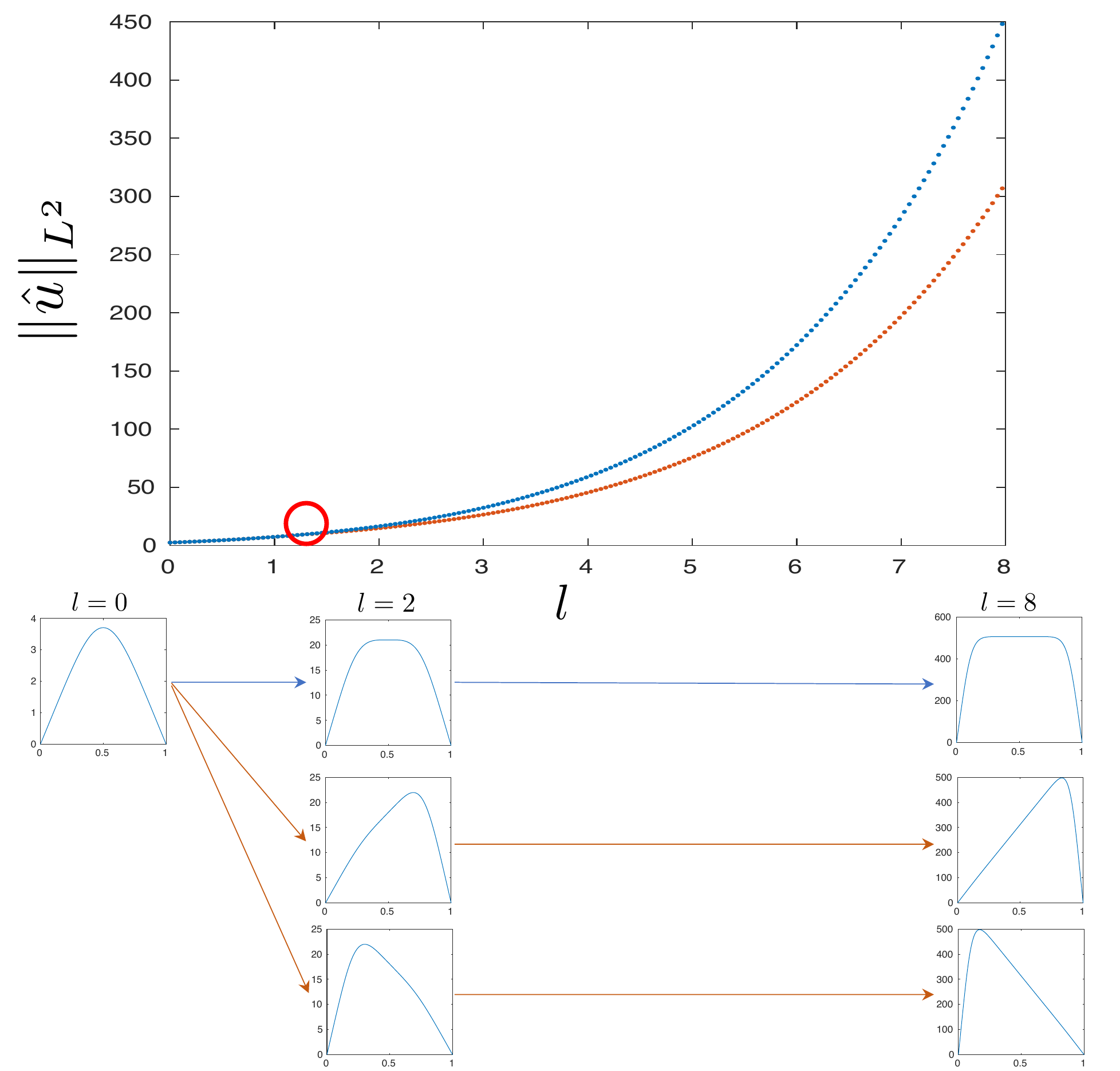}
		\caption{Solution curves for \eqref{henon} on the unit line segment $(0,1)$.}
		\label{1D_SolCur}
	\end{figure}

	Our approximate computation obtained Figure \ref{1D_SolCur}, the solution curve of \eqref{henon} for $0\leq l\leq 8$ ($l$ is always a multiple of 0.05).
	The verified points where $l=0,2,4$ lie on the solution curves.
	According to Figure \ref{1D_SolCur}, a bifurcation point is expected to exist around $[1.20 , 1.25]$.
	

	\begin{table}[H]
		\centering
		\caption{Verification results for $l=0, 2, 4$ on $\Omega = (0,1)$.} 
		\vspace{5truemm}
		\footnotesize 
	\begin{tabular}{l|l|ll|ll}
		\hline
		$l$ & \multicolumn{1}{c|}{0} &\multicolumn{2}{c|}{2} & \multicolumn{2}{c}{4} \\
		\hline
		$\hat{u}$ &
		\begin{minipage}{0.13\textwidth}
			\vspace{2truemm}
			\begin{center}
				\includegraphics[width=1\textwidth]{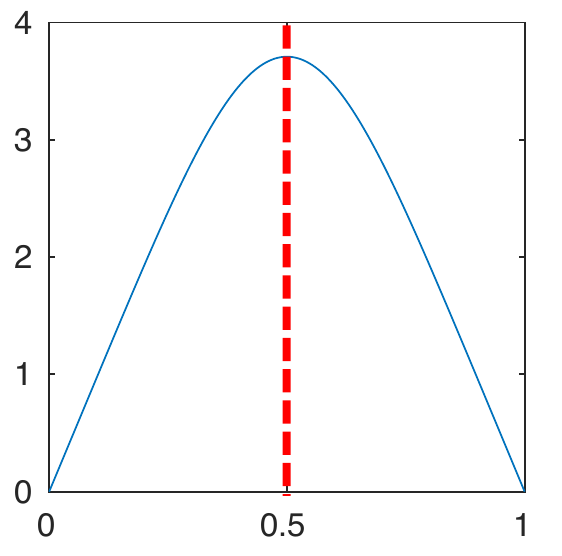}
			\end{center}
			\vspace{-2truemm}
		\end{minipage} &
		
		\begin{minipage}{0.13\textwidth}
			\vspace{2truemm}
			\begin{center}
				\includegraphics[width=1\textwidth]{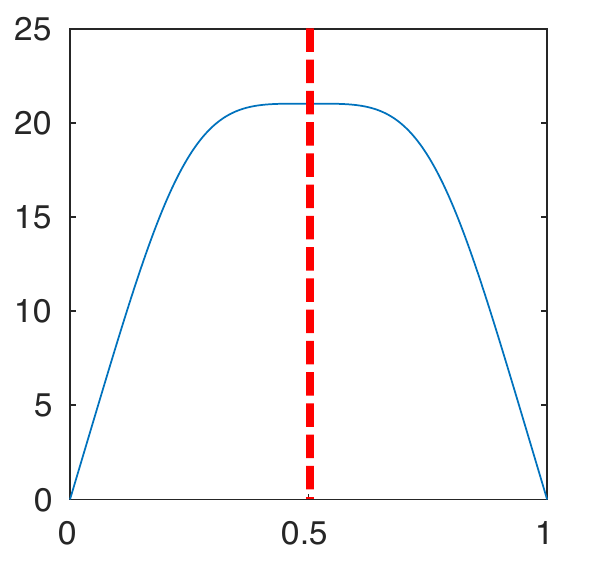}
			\end{center}
			\vspace{-2truemm}
		\end{minipage} &
		
		\begin{minipage}{0.13\textwidth}
			\vspace{2truemm}
			\begin{center}
				\includegraphics[width=1\textwidth]{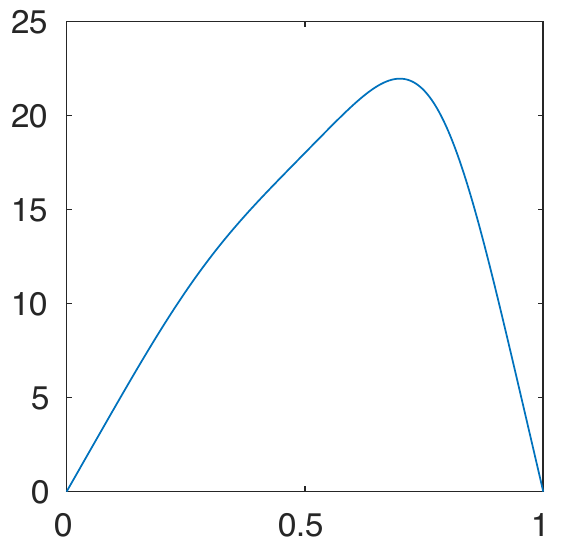}\\
				\includegraphics[width=1\textwidth]{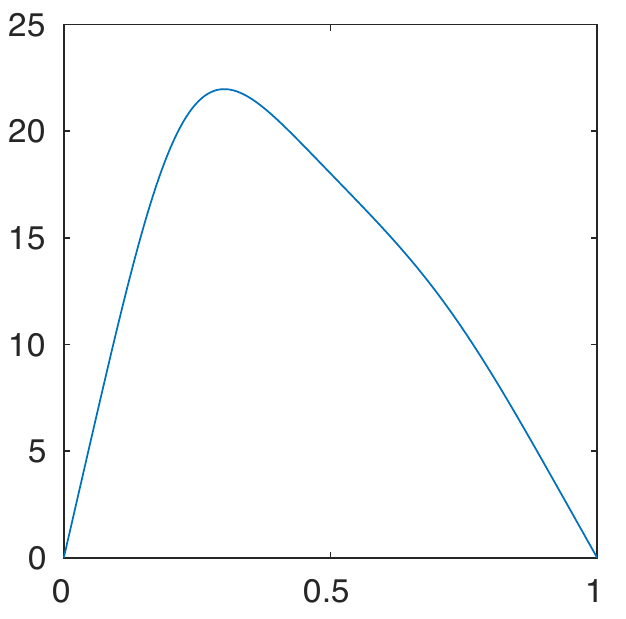}
			\end{center}
			\vspace{-2truemm}
		\end{minipage} &
		
		\begin{minipage}{0.13\textwidth}
			\vspace{2truemm}
			\begin{center}
				\includegraphics[width=1\textwidth]{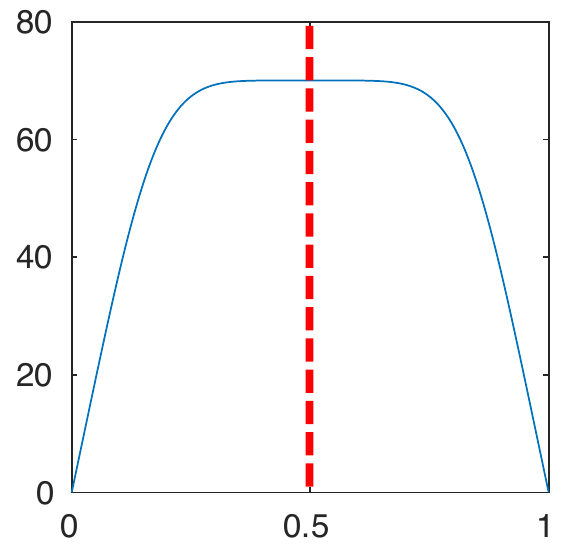}
			\end{center}
			\vspace{-2truemm}
		\end{minipage} &
		
		\begin{minipage}{0.13\textwidth}
			\vspace{2truemm}
			\begin{center}
				\includegraphics[width=1\textwidth]{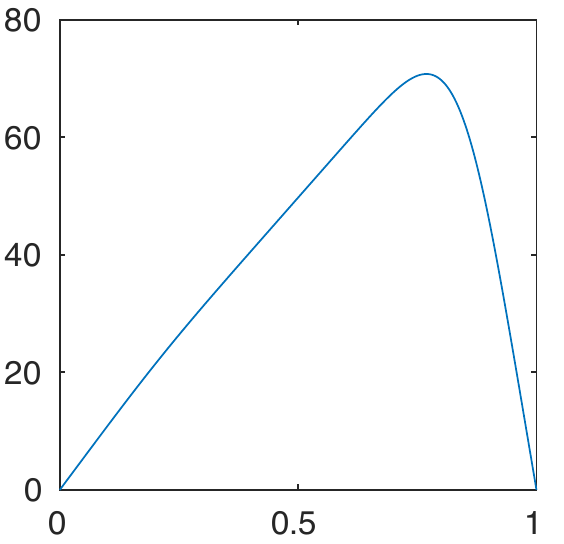}\\
				\includegraphics[width=1\textwidth]{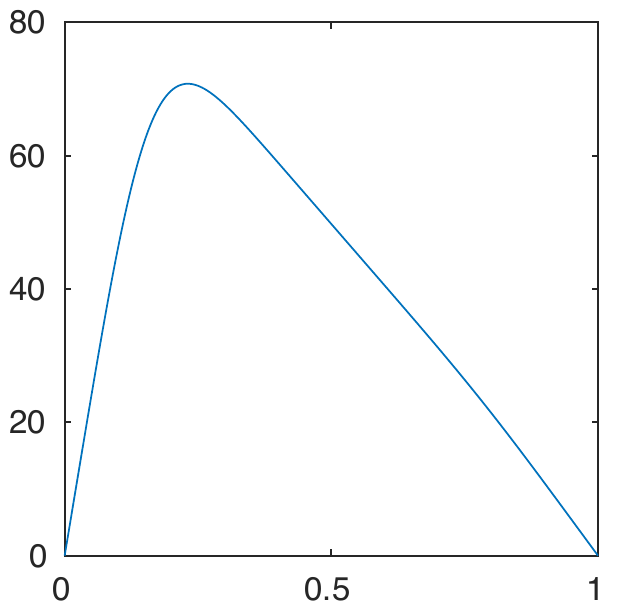}
			\end{center}
			\vspace{-2truemm}
		\end{minipage}
		
		\\ \hline 
		Solution space & $V^1$ & $V^1$ & $V$ & $V^1$ & $V$ \\
		$M_u$ & 40 & 40 & 40 & 40 & 40 \\
		$M$ & 40 & 40 & 40 & 40 & 40 \\
		\footnotesize{$\|F(\hat{u})\|_{H^{-1}}$ } & 2.95468e-12 & 8.35842e-8 & 4.03869e-6 & 9.25374e-6 & 3.36995e-4 \\
		\footnotesize{$\|F'^{-1}_{\hat{u}}\|_{\mathcal{L}(H^{-1},H^1_0)}$} & 2.02207 & 4.19470 & 3.25043 & 1.82276 & 2.16009 \\
		$L$ & 1.28660 & 2.04106 & 1.89034 & 1.78289 & 1.47489 \\
		$\alpha$ & 5.97456e-12 & 3.50610e-7 & 1.31275e-5 & 1.68674e-5 & 7.27937e-4 \\
		$\beta$ & 2.60158 & 8.56162 & 6.14441 & 3.24977 & 3.18587 \\
		$r_A$ & 6.04051e-12 &  4.15274e-7 & 1.51947e-5 & 2.06429e-5 & 9.27220e-4 \\ 
		$r_R$ & 7.60887e-13 & 7.72615e-9 & 2.89288e-7 & 1.00215e-7 & 4.97806e-6 \\ 
		Peak & 3.70815 & 21.0522 & 22.0954 & 70.3607 & 71.2910 \\  \hdashline
		$\mu_{1}$ & -1.99999 & -2.00000 & -1.99999 & -1.99999 & -1.99999 \\
		$\mu_{2}$ & 0.500000 & -0.238397 & 0.356085 & -0.657337 & 0.588997 \\
		$\mu_{3}$ & 0.800001 & 0.703809 & 0.679874 & 0.671403 & 0.755696 \\
		$\mu_{4}$ & 0.892858 & 0.783274 & 0.865471 & 0.733254 & 0.859840 \\
		$\mu_{5}$ & 0.933334 & 0.894429 & 0.910538 & 0.880449 & 0.924964 \\
		\hline
		\end{tabular}
		\label{1D}
	\end{table}
	\footnotesize{
		\noindent
		Solution space: $V:=H^1_0(\Omega)$ and the subspace $V^1$ is defined by \eqref{eq:V1_def_1D}\\
		$M_u$: number of basis functions for constructing approximate solution $\hat{u} \in V_{M_u}$ or $\hat{u} \in V^{1}_{M_u}$ \\
		$M$: number of basis functions for calculating $\lambda^M$\\
		$\|F(\hat{u})\|_{H^{-1}}$: upper bound for the residual norm estimated via \eqref{eq:res}\\
		$\|F'^{-1}_{\hat{u}}\|_{\mathcal{L}(H^{-1},H^1_0)}$: upper bound for the inverse operator norm estimated via Theorem \ref{invtheo}\\
		$L$: upper bound for Lipschitz constant satisfying \eqref{Lip-satis}\\
		$\alpha$: upper bound for $\alpha$ required in Theorem \ref{theo:nk} \\
		$\beta$: upper bound for $\beta$ required in Theorem \ref{theo:nk} \\
		$r_A$: upper bound for absolute error $\| u- \hat{u} \|_{H^1_0}$\\ 
		$r_R$: upper bound for relative error $\| u- \hat{u} \|_{H^1_0}/\| \hat{u} \|_{H^1_0}$ \\ 
		Peak:  upper bound for the maximum values of the corresponding approximation \\
		$\mu_{1}$--$\mu_{5}$:
		approximations of the five smallest eigenvalues of \eqref{eq:eigforTable}
	}
	
	\renewcommand{\arraystretch}{1}
	\normalsize

    \subsection{Numerical results on the unit square}
    We apply our method to $\Omega=(0,1)^2$ in this subsection.
    As in Subsection \ref{sec:1Dresult}, we again restrict solution spaces and their finite-dimensional subspaces to look for symmetric solutions.
	The following sub-solution spaces of $H^1_0(\Omega)$ are endowed with the same topology:
	\begin{align*}
		V^1 &:= \left\{u \in H^1_0(\Omega): u \text{ is symmetric with respect to } x=\frac{1}{2}\right\},\\
		V^2 &:= \left\{u \in H^1_0(\Omega): u \text{ is symmetric with respect to } y=x \right\},\\
		V^3 &:= 
		\left\{u \in H^1_0(\Omega): u \text{ is symmetric with respect to } y=x \text{ and } y=-x+1 \right\},\\
		V^4 &:= 
		\left\{u \in H^1_0(\Omega): u \text{ is symmetric with respect to } x=\frac{1}{2} , y=\frac{1}{2} , y=x \text{, and } y=-x+1 \right\}.
	\end{align*}
	Then, using $\phi_i$ defined in \eqref{eq:Legendrepoly}, we construct finite-dimensional subspaces $V^{i}_{M}~(M\geq 2)$ for each $V^{i}$ ($i=1,2,3,4$) as 
	\begin{align*}
		V^1_{M} &:= \left\{\sum_{i=1\atop \text{$i$ is odd}}^{M} ~\sum_{j=1}^{M} ~u_{i, j} \phi_{i}(x) \phi_{j}(y)~:~ u_{i, j} \in \mathbb{R} \right\},\\
		V^2_{M} &:= \left\{~\sum_{i=1}^{M} ~~\sum_{j=i}^{M} ~u_{i, j} \psi_{i,j}(x,y)~:~  u_{i, j} \in \mathbb{R} \right\},\\
		V^3_{M} &:= \left\{\sum_{i=1\atop \text{$i$ is odd}}^{M} \sum_{j=i\atop \text{$j$ is odd}}^{M} u_{i, j} \psi_{i,j}(x,y) + \sum_{i=2\atop \text{$i$ is even}}^{M} \sum_{j=i\atop \text{$j$ is even}}^{M} u_{i, j} \psi_{i,j}(x,y)~:~  u_{i, j} \in \mathbb{R} \right\},\\
		V^4_{M} &:= \left\{\sum_{i=1\atop \text{$i$ is odd}}^{M} \sum_{j=i\atop \text{$j$ is odd}}^{M} u_{i, j} \psi_{i,j}(x,y)~:~ u_{i, j} \in \mathbb{R} \right\},
	\end{align*}
	where $\psi_{i,j}$ is defined as
	\begin{align*}
		\psi_{i,j}(x,y):= \phi_{i}(x) \phi_{j}(y)+ \phi_{j}(x) \phi_{i}(y),~~(x,y) \in \Omega,
	\end{align*}
	which is symmetric with respect to the line $y=x$.
	Note that we use the same notation $V^1$ and $V^1_M$ with different meanings than in Subsection \ref{sec:1Dresult}.
	The method presented in Sections \ref{sec:NV} and \ref{sec:eva} can be directly applied when the function spaces $H^{1}_{0}(\Omega )$ and $V_{M}$ are replaced with $V^i$ and $V^{i}_M$, respectively.
	In the solution space $V^{i}_{M}$, approximate solutions $\hat{u}$ were obtained by solving the matrix equation 
	\begin{align}
		{\rm Find}&~\hat{u} \in V^{i}_{M}~{\rm s.t.}~	
		\left(\nabla \hat{u},\nabla v_{M}\right)_{L^2}
		=\left(f(\hat{u}),v_{M}\right)_{L^2}~~{\rm for~all}~v_{M}\in V^{i}_{M} \label{eq:newtonapp}
	\end{align}
	via the usual Newton method.
	Restricting solution spaces reduces the amount of calculation for the same reasons as described in Subsection \ref{sec:1Dresult}.
	Using the evaluation \eqref{eq:Lip_end} when $\Omega=(0,1)^2$ with the center $\boldsymbol{x_0}=(1/2,1/2)$, we evaluated the Lipschitz constant $L$ as
	\begin{align*}
		L \leq 6 \left( \frac{1}{\sqrt{2}} \right) ^l   C^3_{4}(\| \hat{u} \| _{L^{4}}+C_{4}r).
	\end{align*}

	Tables \ref{l=02} and \ref{l=4} show the approximate solutions together with their verification results.
	The red dashed lines indicate the symmetry of each solution.
	We again set $\tau\approx 4.9407\times 10^{-324}$, the minimal positive floating-point number after zero.
	In the tables, $\|F(\hat{u})\|_{H^{-1}}$, $\|F'^{-1}_{\hat{u}}\|_{\mathcal{L}(H^{-1},H^1_0)}$, $L$, $\alpha$, and $\beta$ denote the constants required by Theorem \ref{theo:nk}.
	Moreover, $r_A$ and $r_R$ denote an upper bound for absolute error $\| u- \hat{u} \|_{H^1_0}$ and relative error $\| u- \hat{u} \|_{H^1_0}/\| \hat{u} \|_{H^1_0}$, respectively.
	The values in row ``Peak'' represent upper bounds for the maximum values of the corresponding approximate solutions. 
	We see that error bounds are affected by the number of peaks --- fewer peaks lead to larger error bounds.
	As $l$ increases, the peaks approach the corners of the domain and become higher.
	Therefore, a larger $l$ makes verification based on Theorem \ref{theo:nk} more difficult.
	We succeeded in proving the existence of solutions in all cases in which $l=0,2,4$, including three-peak solutions not found in \cite{3}.
	
	The values in rows $\mu_{1}$--$\mu_{5}$ represent approximations of the five smallest eigenvalues of \eqref{eq:eigforTable} discretized in $V_{30} \subset H^1_0(\Omega)$, which is spanned by the basis functions $\phi_n$ ($n=1,2,\cdots,30$) without the restriction of symmetry.
	When $l=4$, the number of negative eigenvalues $\mu$ coincides with the number of peaks.
	
\begin{figure}[H]
		\centering
		\includegraphics[width=0.8\textwidth]{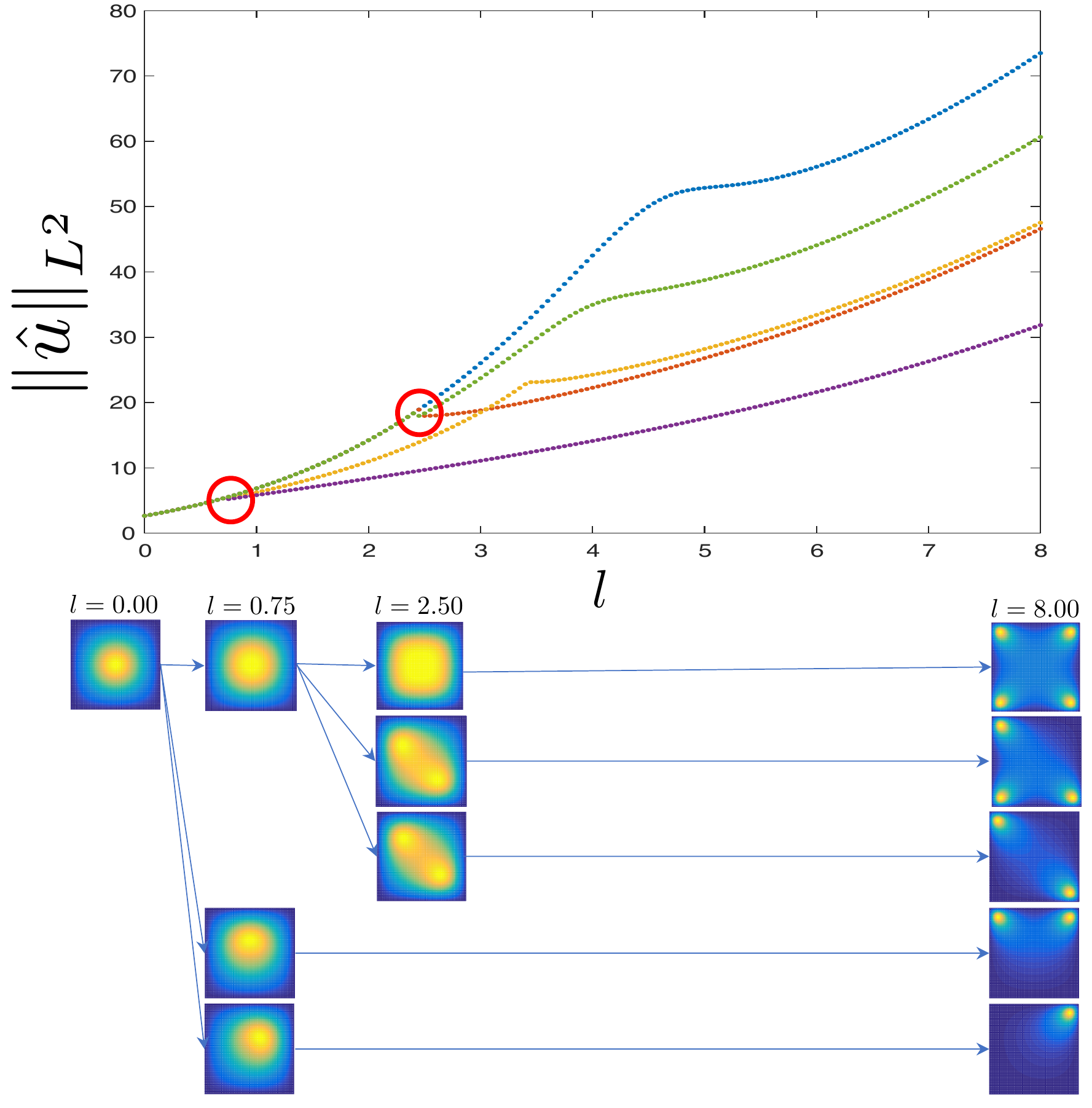}
		\caption{Solution curves for \eqref{henon} on the unit square $(0,1)^2$.}
		\label{zu2}
	\end{figure}

	Our approximate computation obtained Figure \ref{zu2}, the solution curves of \eqref{henon} for $0\leq l\leq 8$ ($l$ is always a multiple of 0.05).
	If the vertical axis scaling is changed,
	the curves coincide with those in \cite[Figure 2]{3} except for that corresponding to the three-peak solutions after the point around $[2.35 , 2.40]$.
	The verified points where $l=0,2,4$ lie on the solution curves.
	According to Figure \ref{zu2}, two bifurcation points are expected to exist around $[0.55 , 0.60]$ and $[2.35 , 2.40]$.
	We expect the single-solution curve bifurcates to three at the first bifurcation point around $[0.55 , 0.60]$, and then one of them further bifurcates to three at the second point around $[2.35 , 2.40]$.
	
	\renewcommand{\arraystretch}{1.2}
	\begin{table}[H]
		\centering
		\caption{Verification results  for $l=0,2$ on the unit square $(0,1)^2$.} 
		\vspace{5truemm}
		\footnotesize
		\begin{tabular}{l|l|lll}
			\hline
			$l$ & \multicolumn{1}{c|}{0} & \multicolumn{3}{c}{2}   \\
			\hline
			3D $\hat{u}$ &
			\begin{minipage}{0.19\textwidth}
				\centering
				\includegraphics[width=0.825\textwidth, bb=110 275 470 560,clip]{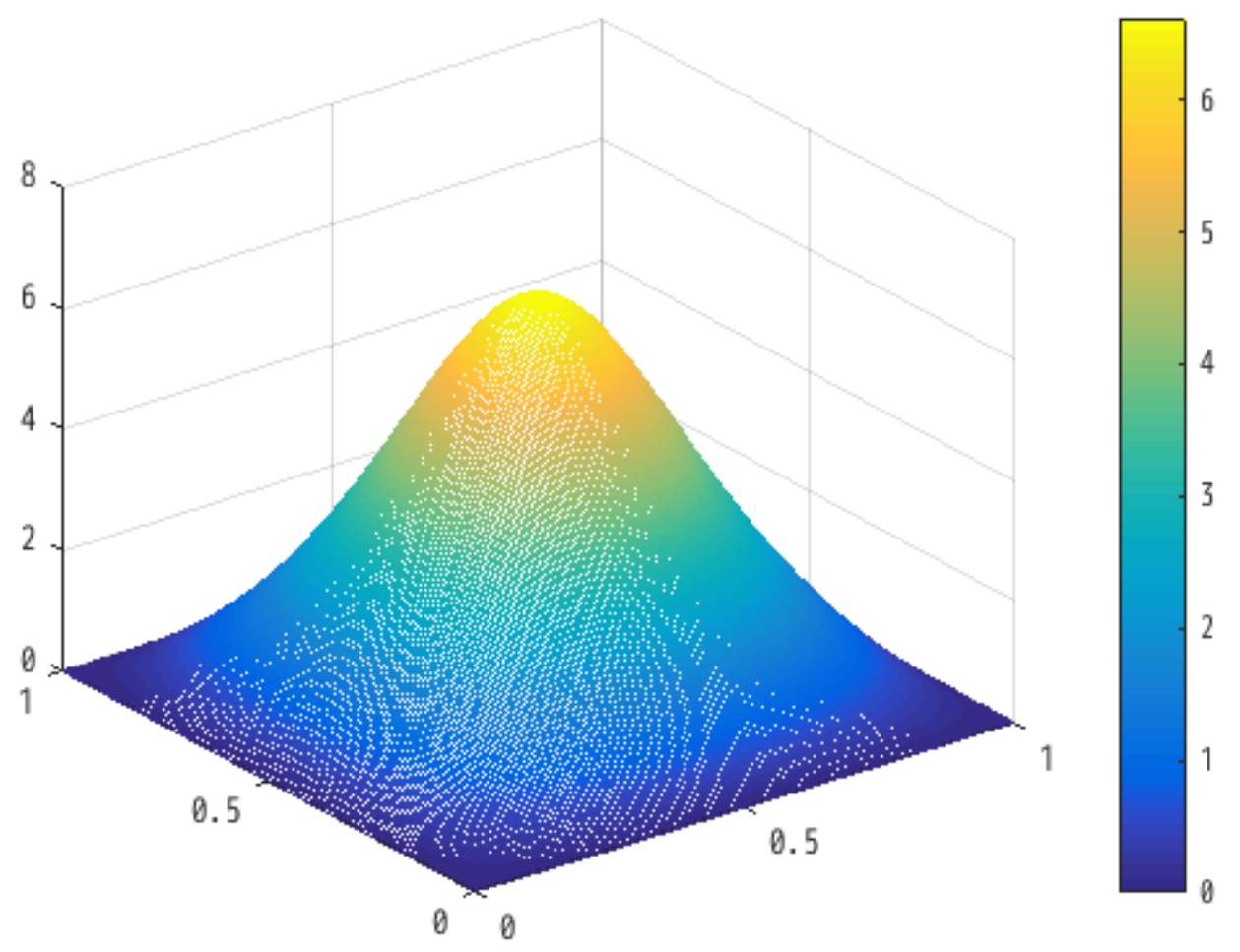}
			\end{minipage} &
			
			\begin{minipage}{0.19\textwidth}
				\centering
				\includegraphics[width=0.825\textwidth, bb=110 275 470 560,clip]{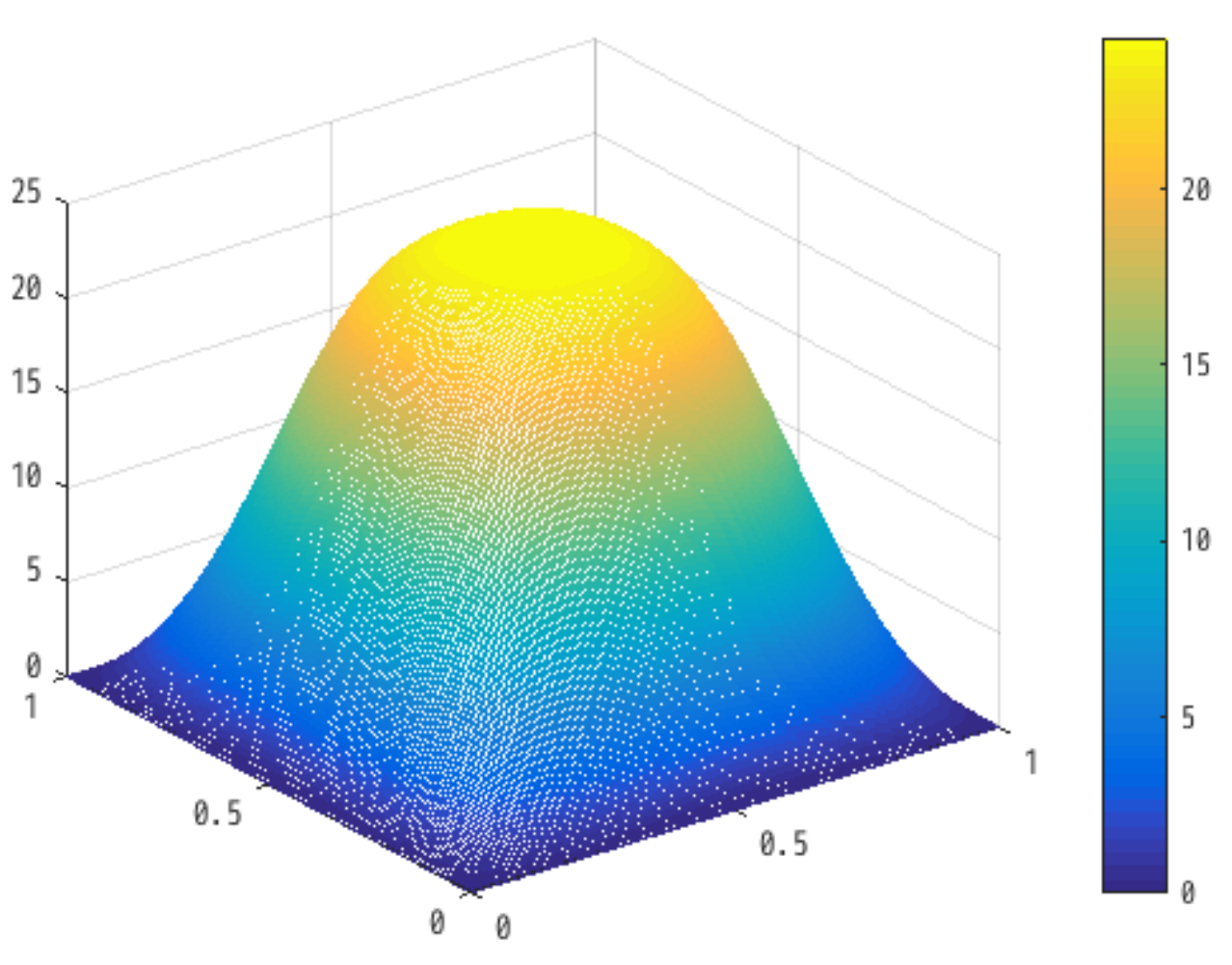}
			\end{minipage} &
			
			\begin{minipage}{0.19\textwidth}
				\centering
				\includegraphics[width=0.825\textwidth, bb=110 275 470 560,clip]{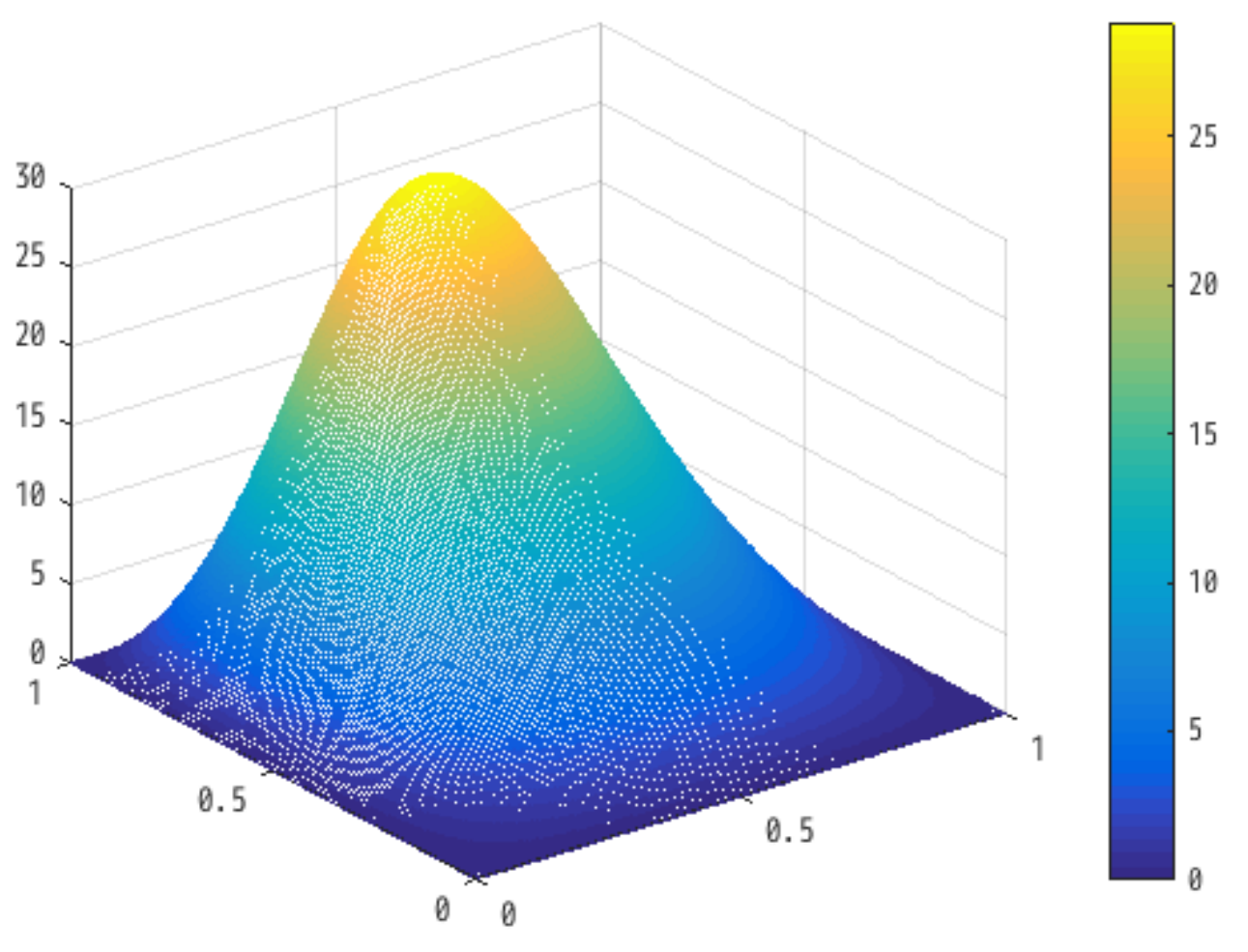}
			\end{minipage} &
			
			\begin{minipage}{0.14\textwidth}
				\centering
				\includegraphics[width=1.1\textwidth, bb=110 275 470 560,clip]{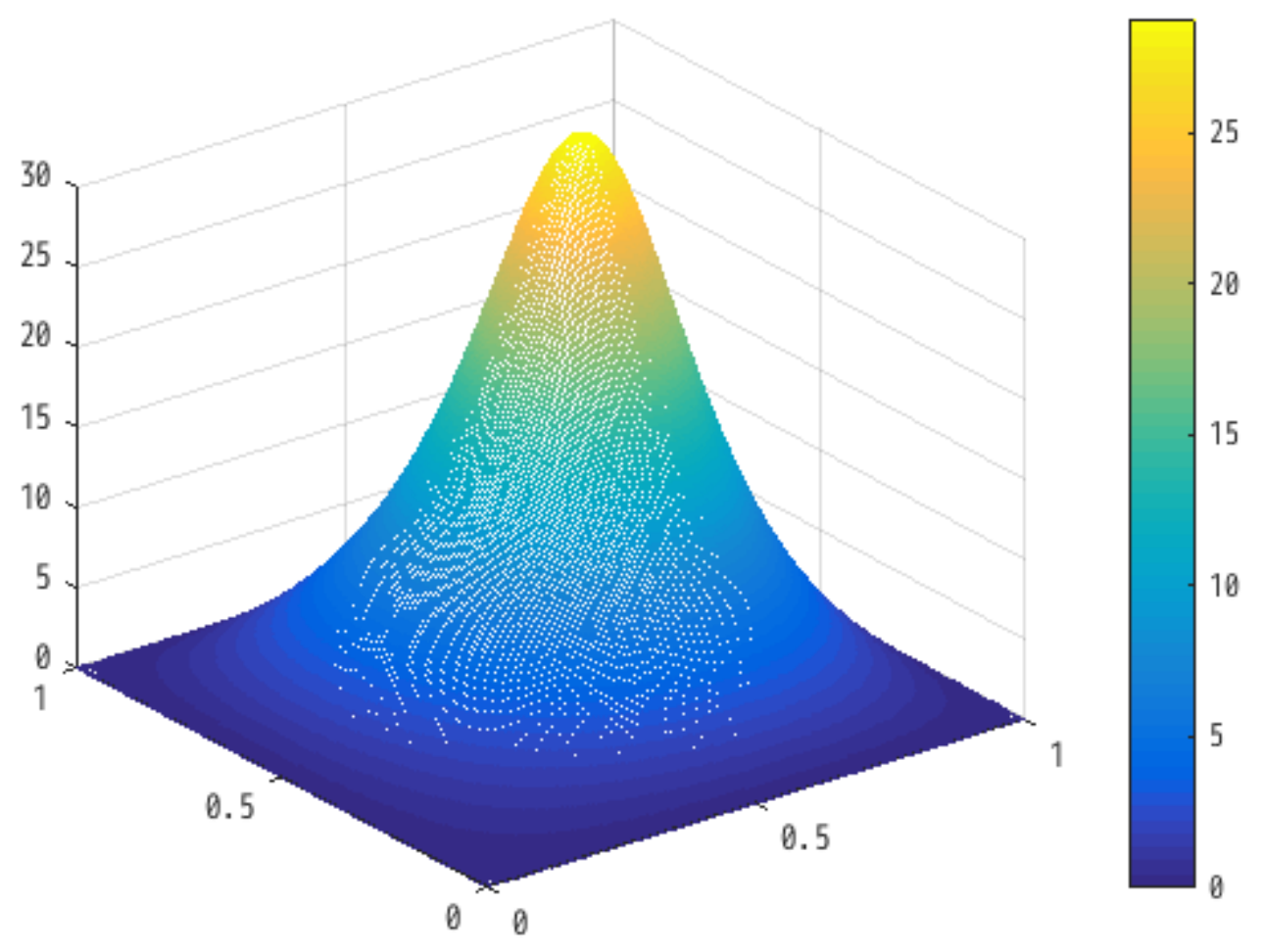}
			\end{minipage} 
			\\	 
			2D $\hat{u}$ &
			~~
			\begin{minipage}{0.14\textwidth}
				\vspace{1truemm}
				\begin{center}
					\includegraphics[width=1\textwidth]{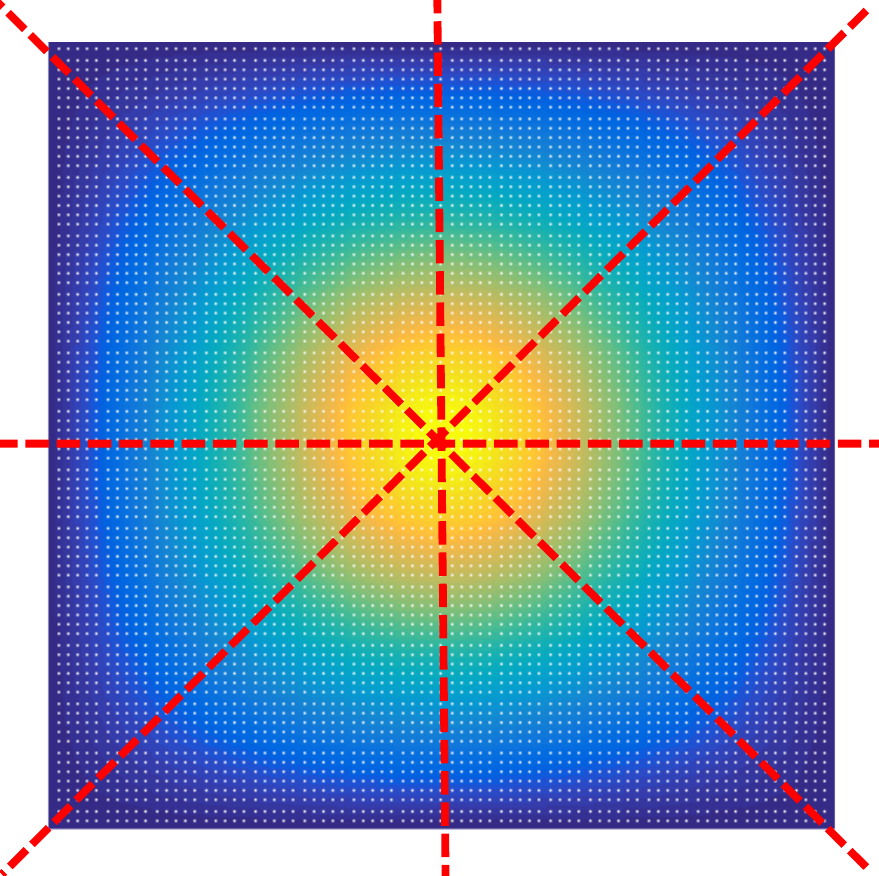}
				\end{center}
				\vspace{-2mm}
			\end{minipage} 
			&
			\multicolumn{1}{c}{
				\begin{minipage}{0.14\textwidth}
					\vspace{1truemm}
					\begin{center}
						\includegraphics[width=1\textwidth]{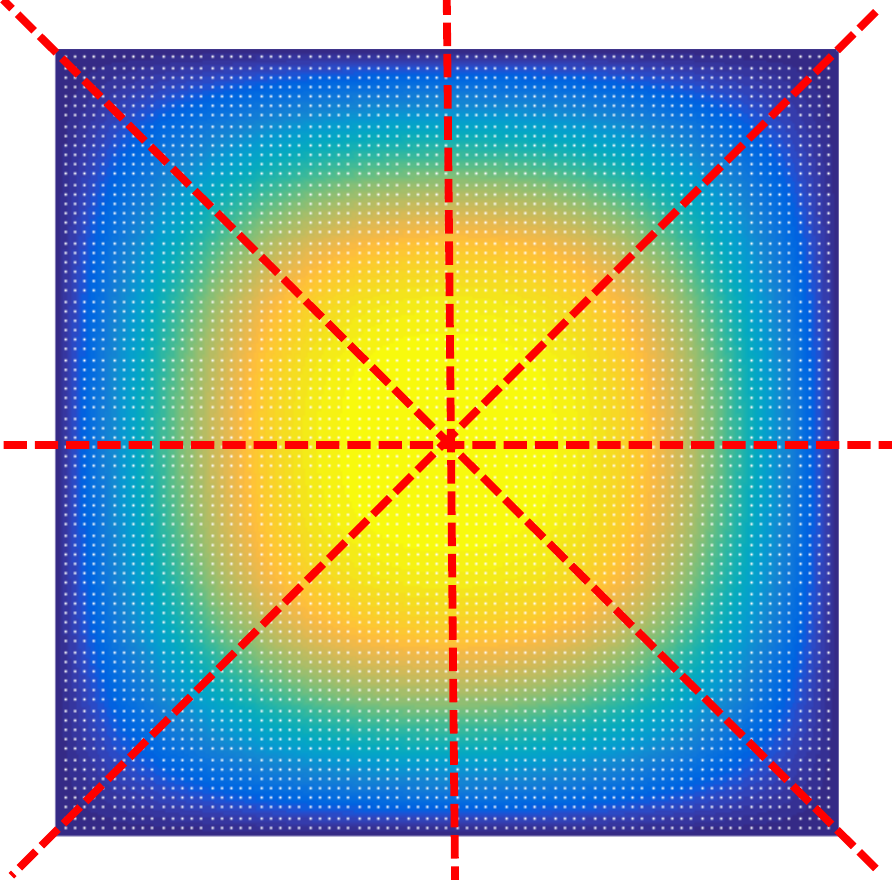}
					\end{center}
					\vspace{-2truemm}
				\end{minipage} 
			}
			&
			\multicolumn{1}{c}{
				\begin{minipage}{0.14\textwidth}
					\vspace{1truemm}
					\begin{center}
						\includegraphics[width=1\textwidth]{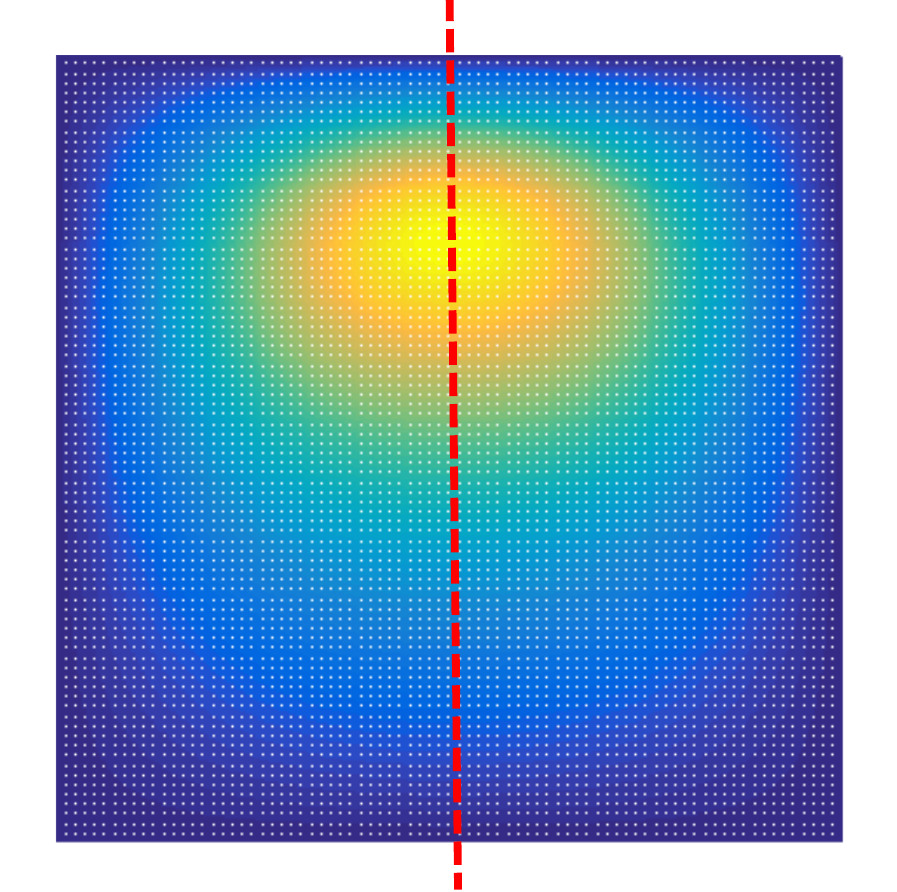}
					\end{center}
					\vspace{-2truemm}
				\end{minipage} 
			}
			&
			\begin{minipage}{0.14\textwidth}
				\vspace{1truemm}
				\begin{center}
					\includegraphics[width=1\textwidth]{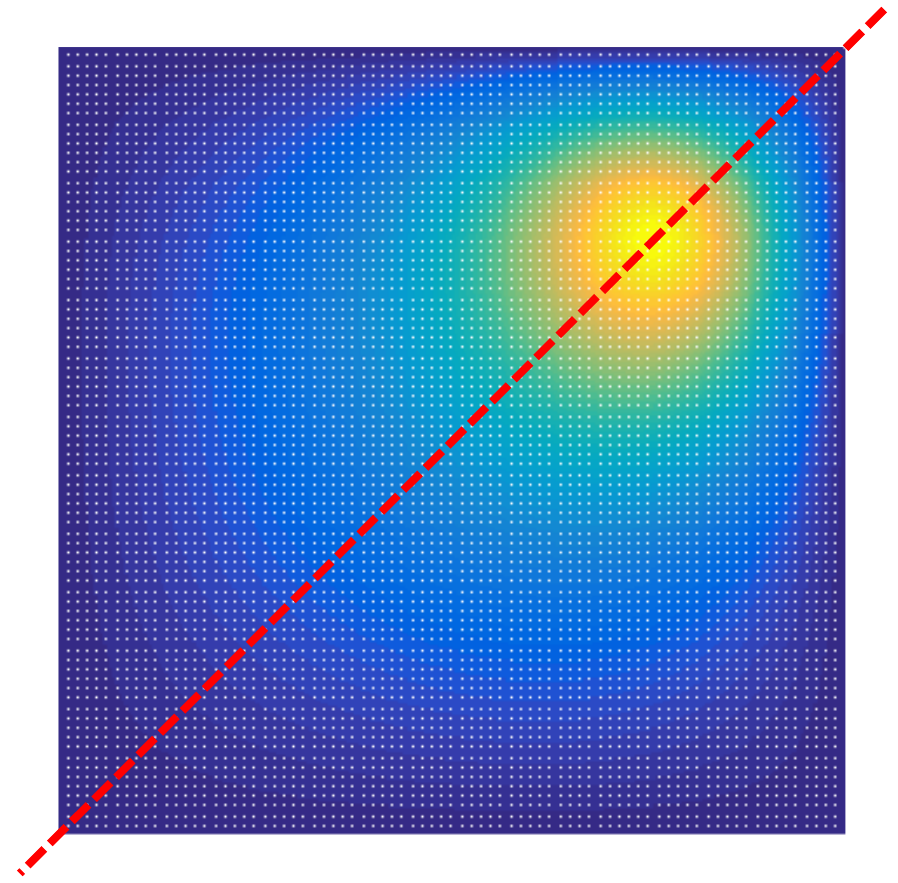}
				\end{center}
				\vspace{-2truemm}
			\end{minipage} 
			\\ \hline 
			Solution space & $V^4$ & $V^4$ & $V^1$ & $V^2$ \\
			$M_u$ & 40 & 40 & 60& 60 \\
			$M$ & 40 & 40 & 40 & 40\\
			\footnotesize{$\|F(\hat{u})\|_{H^{-1}}$ }& 1.17370e-7 & 3.96407e-7 & 1.19312e-8 & 4.22257e-7\\
			\footnotesize{$\|F'^{-1}_{\hat{u}}\|_{\mathcal{L}(H^{-1},H^1_0)}$} & 1.70326 & 2.26200 & 15.19763 & 36.47472\\
			$L$ & 6.78398e-1 & 1.64252 & 1.43209 & 1.21150 \\
			$\alpha$ & 1.99910e-7 & 8.96672e-7 & 1.81325e-7 & 1.54017e-5 \\
			$\beta$ & 1.15549 & 3.71537 & 21.76424 & 44.18887 \\
			$r_A$ & 4.63296e-8 & 2.55597e-7 & 1.44557e-7 & 2.48634e-5  \\ 
			$r_R$ & 3.76958e-9 & 3.98528e-9 &  2.45351e-9 & 4.63166e-7  \\ 
			Peak & 6.62326 & 24.36528 & 29.03437 & 29.20268  \\ \hdashline 
			$\mu_{1}$ & -1.99999 & -1.99999 & -1.99999 & -1.99999 \\
			$\mu_{2}$ & 0.220034 & -0.410090 & -0.273589 & 0.196622 \\
			$\mu_{3}$ & 0.220034 & -0.410090 & 0.233061 & 0.208937 \\
			$\mu_{4}$ & 0.604521 & 0.114826 & 0.457439 & 0.585268 \\
			$\mu_{5}$ & 0.658421 & 0.298974 & 0.517021 & 0.639470 \\
			\hline
		\end{tabular}
		\label{l=02} 
	\end{table}
	\footnotesize{
		\noindent
		Solution space: restricted solution space $V^i \subset H^1_0(\Omega)$ \\
		$M_u$: number of basis functions with respect to $x$ and $y$ for constructing approximate solution $\hat{u} \in V^i_{M_u}$\\
		$M$: number of basis functions with respect to $x$ and $y$ for calculating $\lambda^M$\\
		$\|F(\hat{u})\|_{H^{-1}}$: upper bound for the residual norm estimated via \eqref{eq:res}\\
		$\|F'^{-1}_{\hat{u}}\|_{\mathcal{L}(H^{-1},H^1_0)}$: upper bound for the inverse operator norm estimated via Theorem \ref{invtheo}\\
		$L$: upper bound for Lipschitz constant satisfying \eqref{Lip-satis}\\
		$\alpha$: upper bound for $\alpha$ required in Theorem \ref{theo:nk} \\
		$\beta$: upper bound for $\beta$ required in Theorem \ref{theo:nk} \\
		$r_A$: upper bound for absolute error $\| u- \hat{u} \|_{H^1_0}$\\ 
		$r_R$: upper bound for relative error $\| u- \hat{u} \|_{H^1_0}/\| \hat{u} \|_{H^1_0}$ \\ 
		Peak:  upper bound for the maximum values of the corresponding approximation \\
		$\mu_{1}$--$\mu_{5}$:
		approximations of the five smallest eigenvalues of \eqref{eq:eigforTable}
	}
	\begin{table}[H]
		\centering
		\caption{Verification results  for $l=4$ on the unit square $(0,1)^2$.} 
		\vspace{5truemm}
		\footnotesize 
		\begin{tabular}{l|lllll}
			\hline
			$l$ & \multicolumn{5}{c}{4}  \\
			\hline
			3D $\hat{u}$ &
			\begin{minipage}{0.14\textwidth}
				\centering
				\includegraphics[width=1.1\textwidth, bb=110 275 470 560,clip]{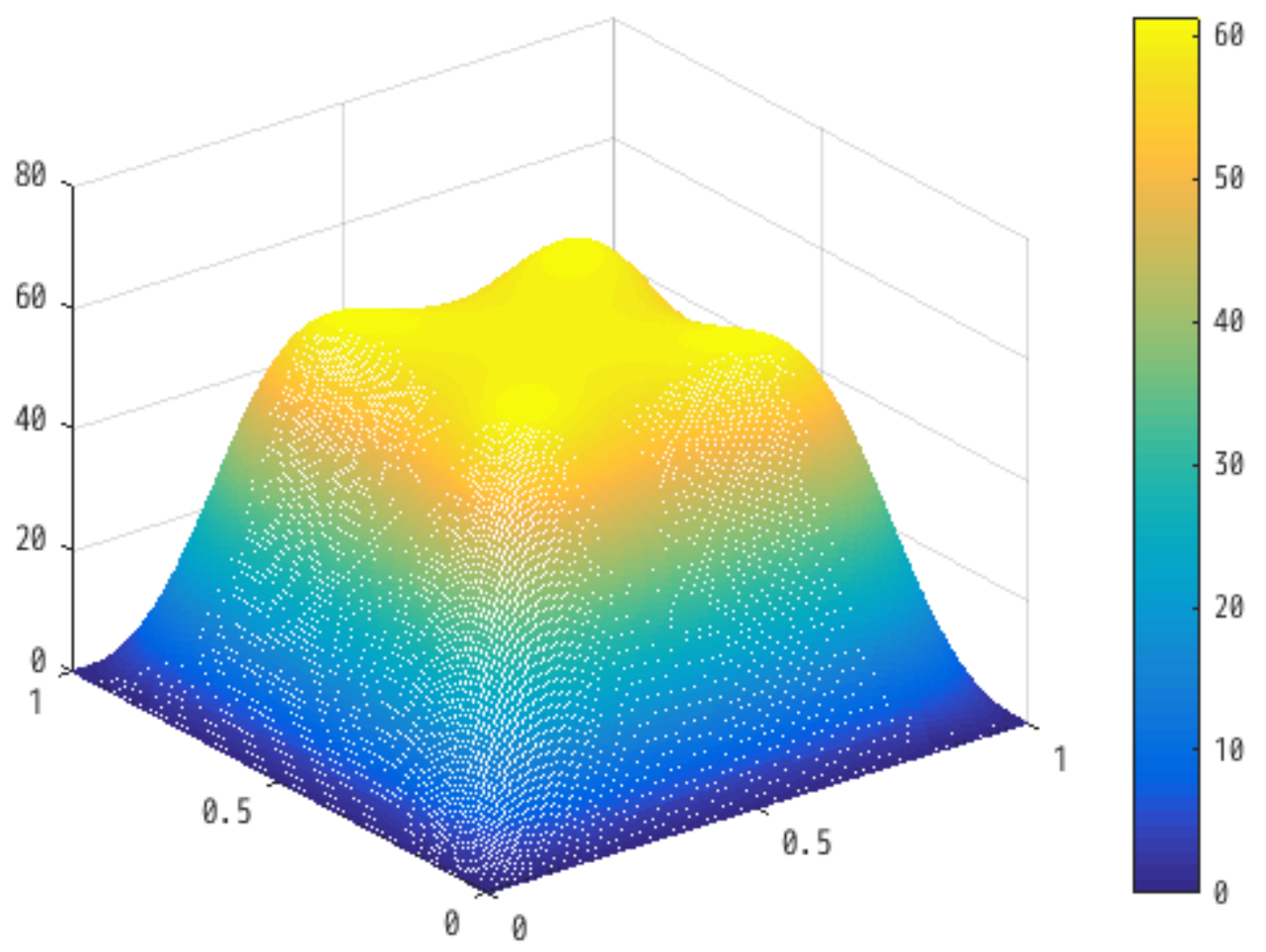}
			\end{minipage} &
			
			\begin{minipage}{0.14\textwidth}
				\centering
				\includegraphics[width=1.1\textwidth, bb=110 275 470 560,clip]{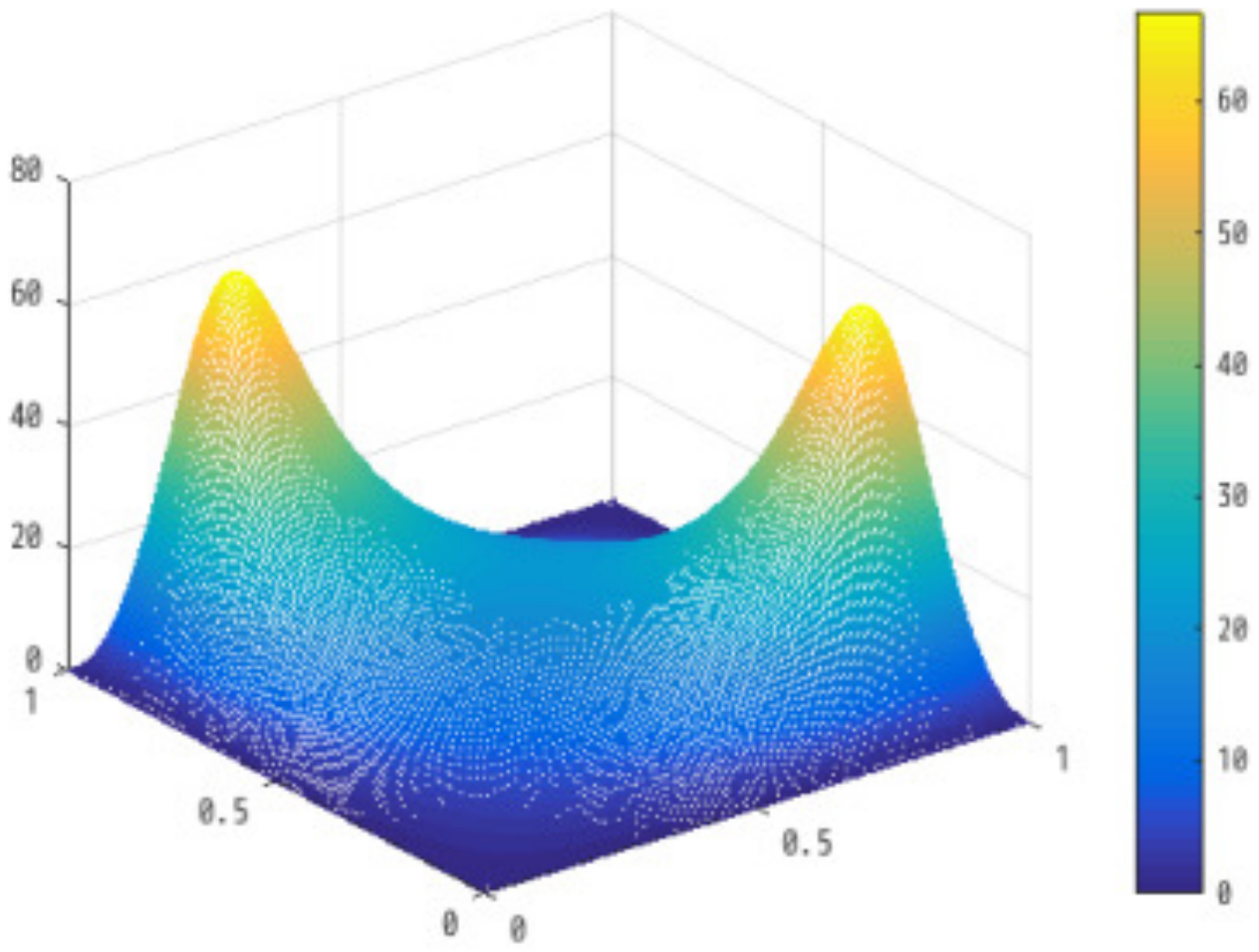}
			\end{minipage} &
			
			\begin{minipage}{0.14\textwidth}
				\centering
				\includegraphics[width=1.1\textwidth, bb=110 275 470 560,clip]{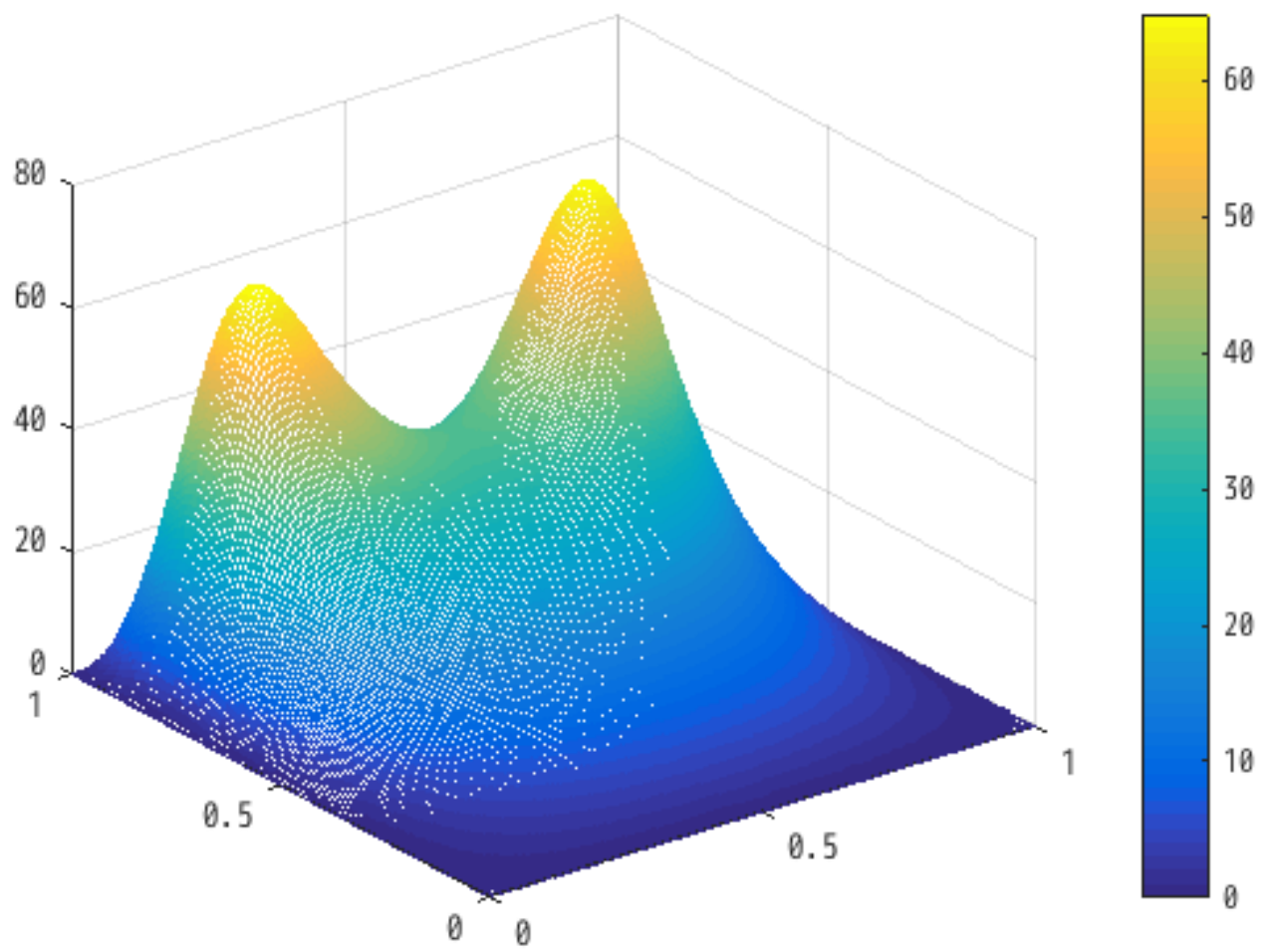}
			\end{minipage} &
			
			\begin{minipage}{0.14\textwidth}
				\centering
				\includegraphics[width=1.1\textwidth, bb=110 275 470 560,clip]{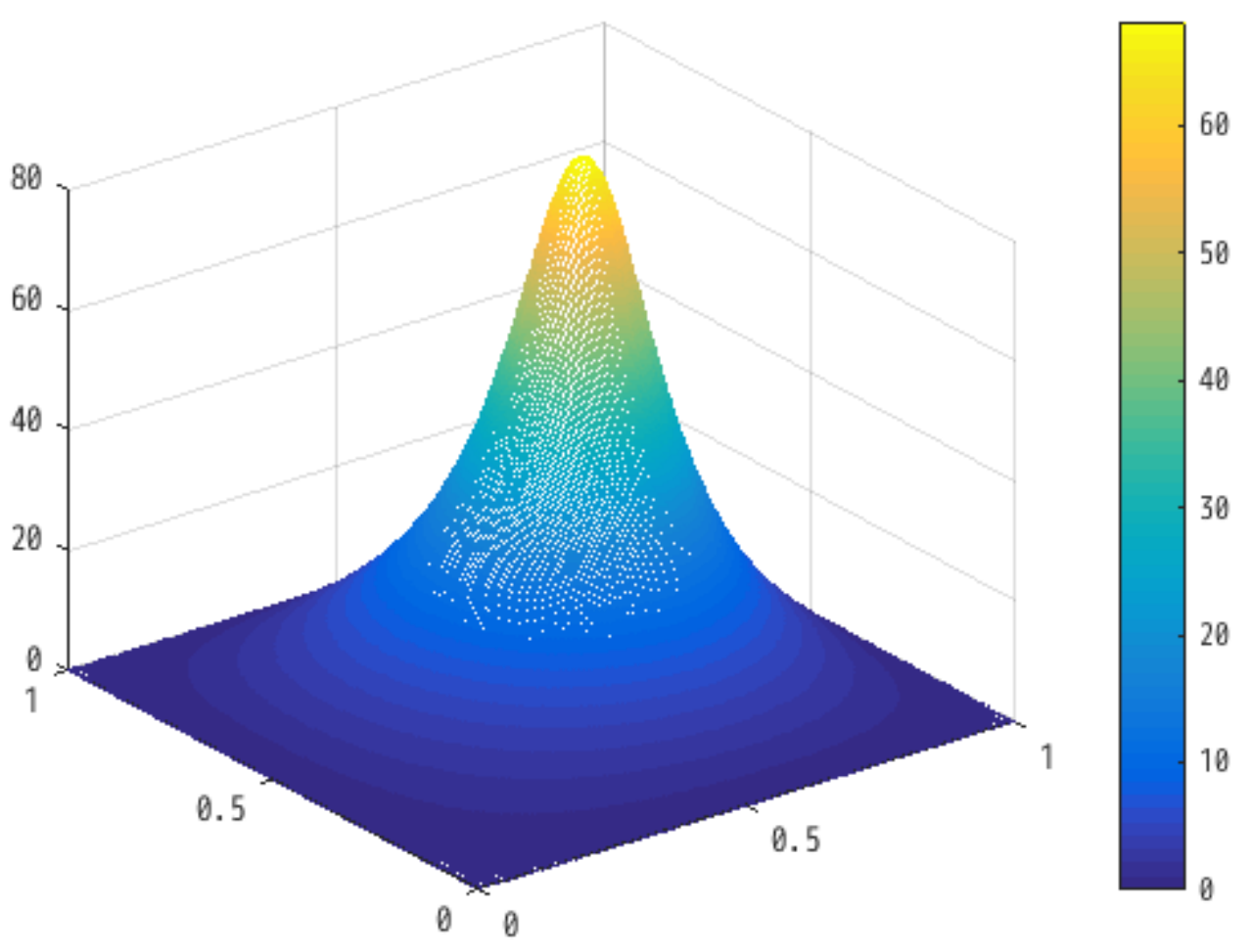}
			\end{minipage} &
			
			\begin{minipage}{0.14\textwidth}
				\centering
				\includegraphics[width=1.1\textwidth, bb=110 275 470 560,clip]{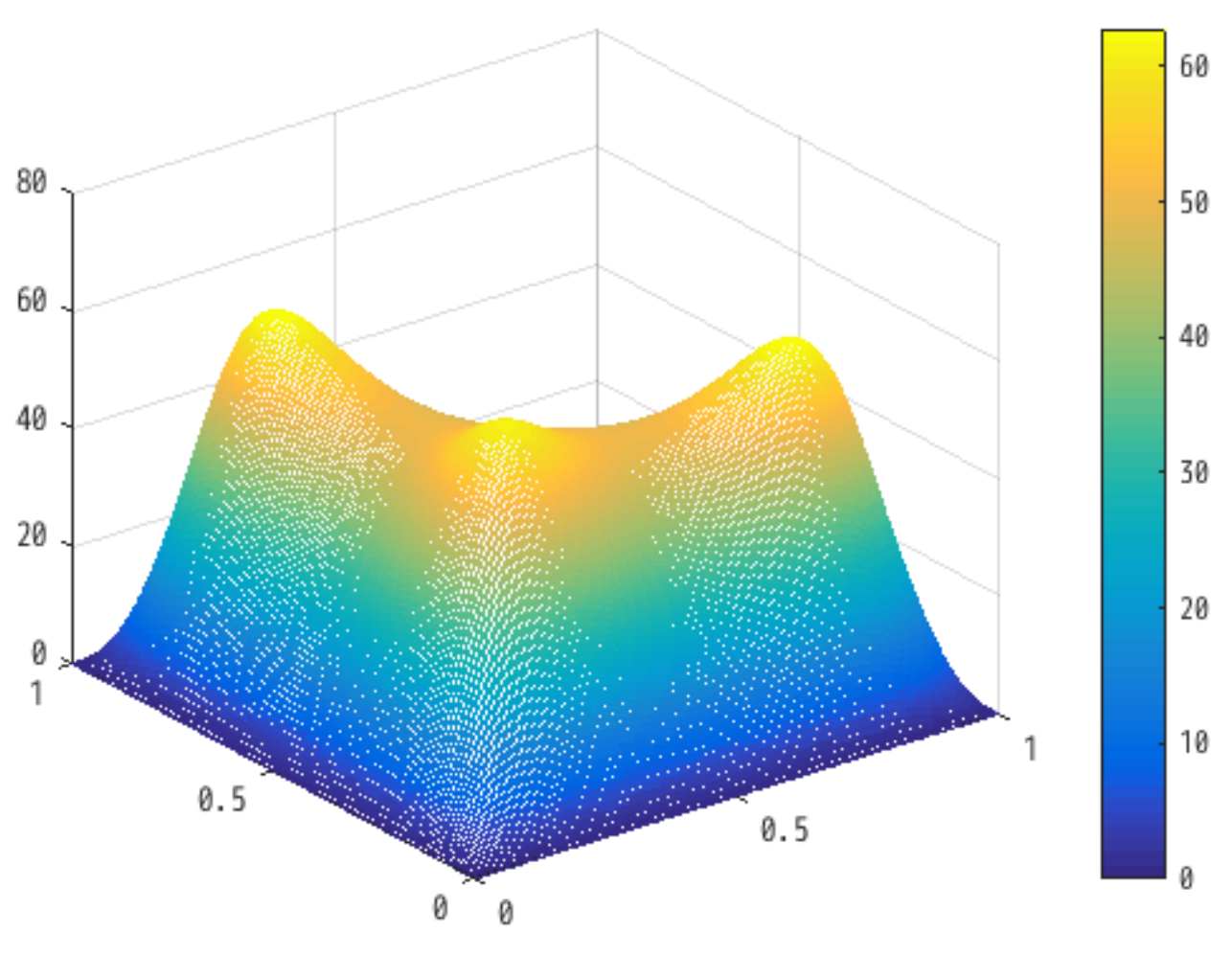}
			\end{minipage}
			
			\\	 
			2D $\hat{u}$ &
			\begin{minipage}{0.14\textwidth}
				\vspace{1truemm}
				\begin{center}
					\includegraphics[width=1\textwidth]{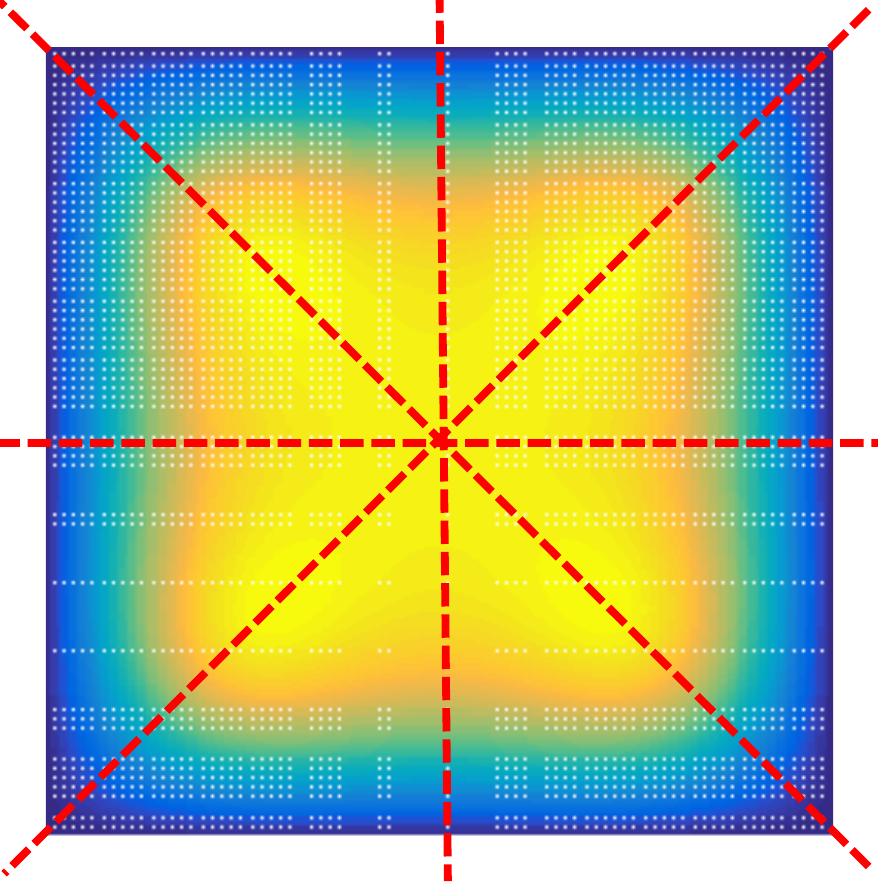}
				\end{center}
				\vspace{-2truemm}
			\end{minipage} &
			
			\begin{minipage}{0.14\textwidth}
				\vspace{1truemm}
				\begin{center}
					\includegraphics[width=1\textwidth]{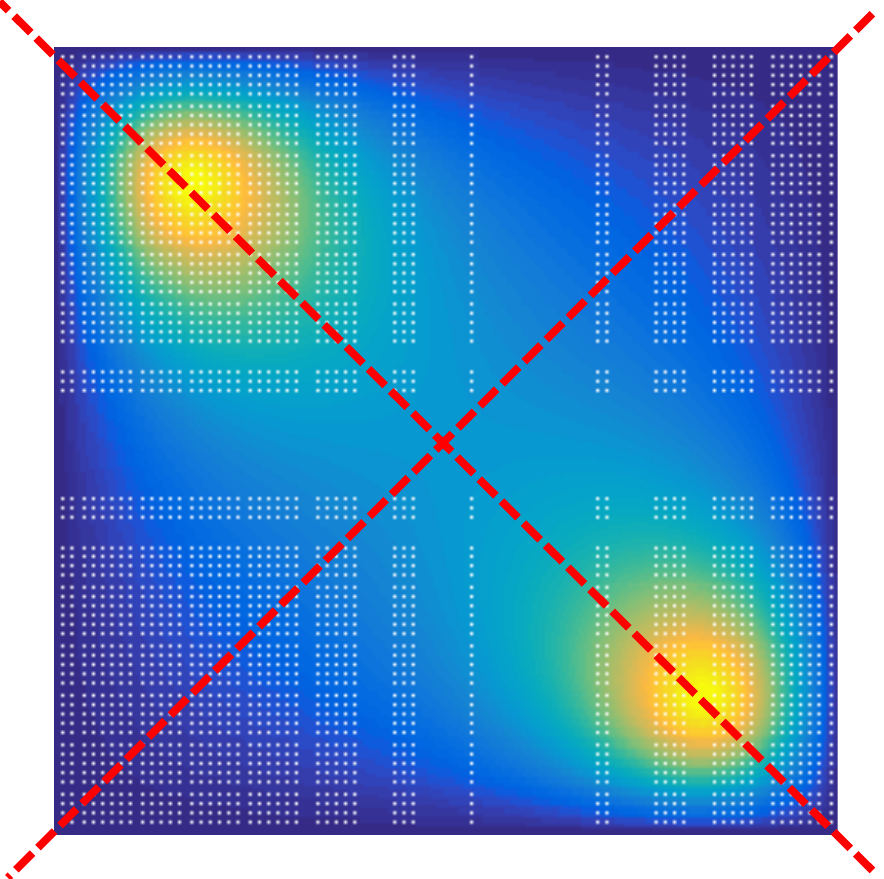}
				\end{center}
				\vspace{-2truemm}
			\end{minipage} &
			
			\begin{minipage}{0.14\textwidth}
				\vspace{1truemm}
				\begin{center}
					\includegraphics[width=1\textwidth]{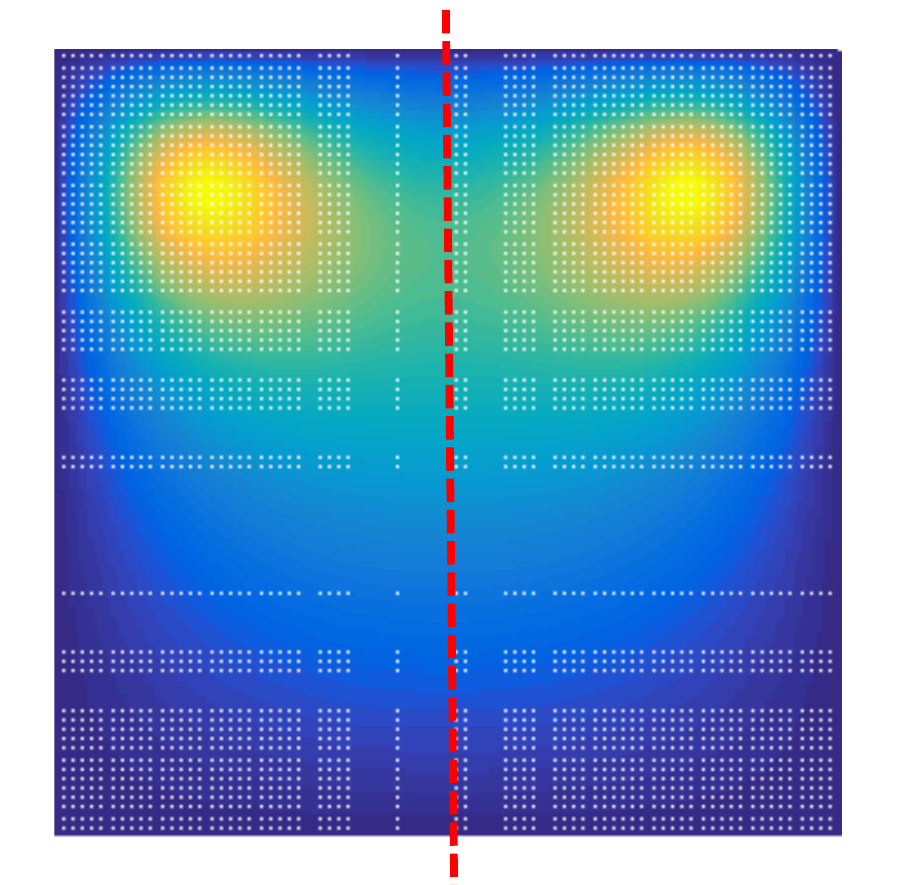}
				\end{center}
				\vspace{-2truemm}
			\end{minipage} &
			
			\begin{minipage}{0.14\textwidth}
				\vspace{1truemm}
				\begin{center}
					\includegraphics[width=1\textwidth]{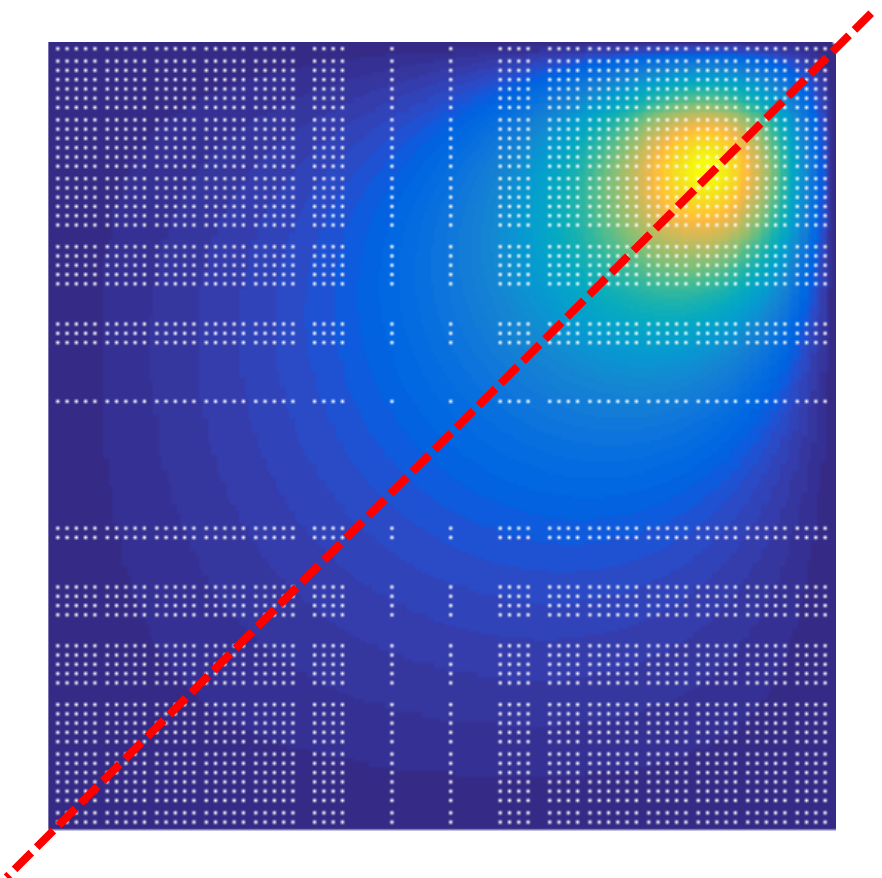}
				\end{center}
				\vspace{-2truemm}
			\end{minipage} &
			
			\begin{minipage}{0.14\textwidth}
				\vspace{1truemm}
				\begin{center}
					\includegraphics[width=1\textwidth]{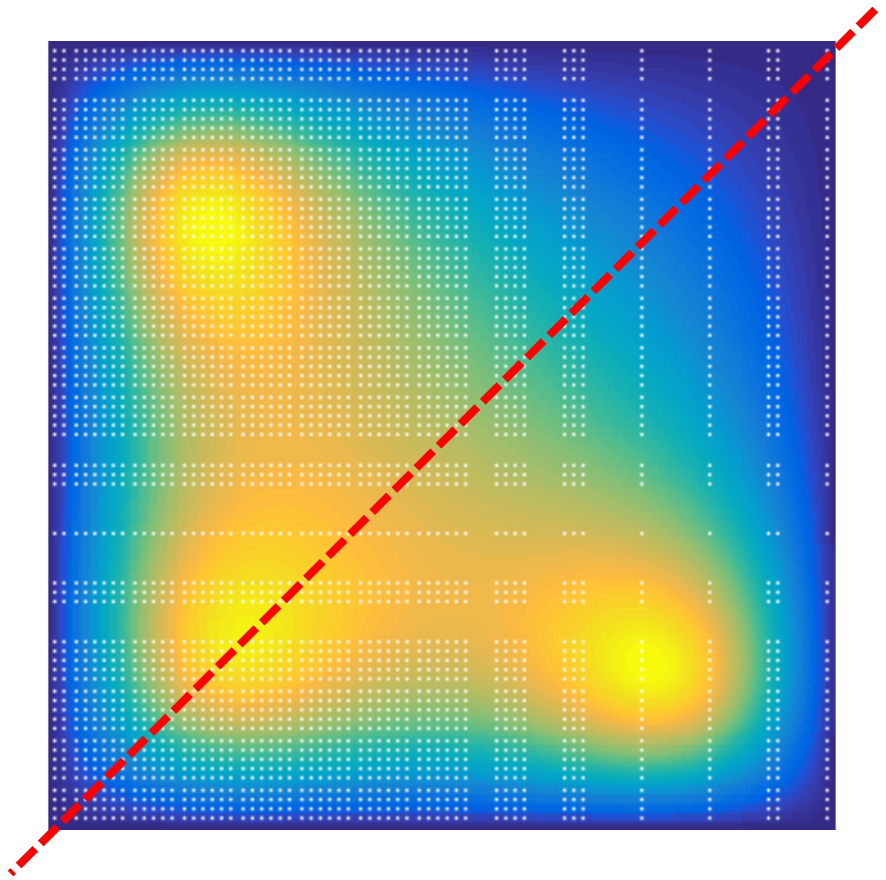}
				\end{center}
				\vspace{-2truemm}
			\end{minipage}
			
			\\ \hline 
			Solution space & $V^4$ & $V^3$ & $V^1$ & $V^2$ & $V^2$ \\
			$M_u$ & 70 & 70 & 70 & 70 & 70 \\
			$M$ & 80 & 80 & 80 & 80 & 80 \\
			\footnotesize{$\|F(\hat{u})\|_{H^{-1}}$ } & 1.88534e-11 & 7.91070e-6 & 4.76970e-7 & 8.47044e-6 & 3.47384e-8 \\
			\footnotesize{$\|F'^{-1}_{\hat{u}}\|_{\mathcal{L}(H^{-1},H^1_0)}$} & 6.82420 & 24.18779 & 78.96665 & 21.26750 & 47.44875 \\
			$L$ & 2.31308 & 1.46531 & 1.55126 & 1.18832 & 1.97091 \\
			$\alpha$ & 1.28659e-10 & 1.91343e-4 & 3.76648e-5 & 1.80145e-4 & 1.64830e-6 \\
			$\beta$ & 15.78486 & 35.44250 & 1.22498e+2 & 25.27251 & 93.51720 \\
			$r_A$ & 4.95952e-11 &  1.73351e-4 & 8.76586e-5 & 1.53306e-4 & 2.32064e-6 \\ 
			$r_R$ & 2.35369e-13 & 9.86681e-7 & 5.12219e-7 & 1.20925e-6 & 1.16657e-8 \\ 
			Peak & 62.30489 & 68.15045 & 66.28947 & 69.69524 & 64.16408 \\ \hdashline
			$\mu_{1}$ & -1.99999 & -1.99996 & -1.99999 & -1.99999 & -1.99999 \\
			$\mu_{2}$ & -0.995156 & -1.86714 & -1.64594 & 0.177691 & -1.46267 \\
			$\mu_{3}$ & -0.995156 & 0.166245 & 0.130875 & 0.251043 & -1.14006 \\
			$\mu_{4}$ & -0.689431 & 0.205039 & 0.253364 & 0.591950 & 0.131828 \\
			$\mu_{5}$ & 0.210478 & 0.258004 & 0.272595 & 0.658008 & 0.175494 \\
			\hline
		\end{tabular}
		\label{l=4}
	\end{table}
	\footnotesize{
		\noindent
		Solution space: restricted solution space $V^i \subset H^1_0(\Omega)$ \\
		$M_u$: number of basis functions with respect to $x$ and $y$ for constructing approximate solution $\hat{u} \in V^i_{M_u}$\\
		$M$: number of basis functions with respect to $x$ and $y$ for calculating $\lambda^M$\\
		$\|F(\hat{u})\|_{H^{-1}}$: upper bound for the residual norm estimated via \eqref{eq:res}\\
		$\|F'^{-1}_{\hat{u}}\|_{\mathcal{L}(H^{-1},H^1_0)}$: upper bound for the inverse operator norm estimated via Theorem \ref{invtheo}\\
		$L$: upper bound for Lipschitz constant satisfying \eqref{Lip-satis}\\
		$\alpha$: upper bound for $\alpha$ required in Theorem \ref{theo:nk} \\
		$\beta$: upper bound for $\beta$ required in Theorem \ref{theo:nk} \\
		$r_A$: upper bound for absolute error $\| u- \hat{u} \|_{H^1_0}$\\ 
		$r_R$: upper bound for relative error $\| u- \hat{u} \|_{H^1_0}/\| \hat{u} \|_{H^1_0}$ \\ 
		Peak:  upper bound for the maximum values of the corresponding approximation \\
		$\mu_{1}$--$\mu_{5}$:
		approximations of the five smallest eigenvalues of \eqref{eq:eigforTable}
	}
	\renewcommand{\arraystretch}{1}
	\section{Conclusion} \label{sec:conclusion}
	\normalsize
	We designed a numerical verification method for proving the existence of solutions of the H\'enon equation \eqref{henon} on a bounded domain based on the Newton-Kantorovich theorem.
	We applied our method to the domains $\Omega=(0,1)^N$ $(N=1,2)$,
	proving the existence of several solutions of \eqref{henon} nearby a numerically computed approximation $\hat{u}$.
	In particular, we found a set of undiscovered solutions with three peaks on the square domain $\Omega=(0,1)^2$.
	Approximate computations generated the solution curves of \eqref{henon} for $0\leq l\leq 8$ in Figures \ref{1D_SolCur} and \ref{zu2}.
	Future work should verify the existence of solutions for arbitrary $l \in [0,a]$, given a large $a>0$, and prove the bifurcation structure for \eqref{henon} in a strict mathematical sense.
	\section{Acknowledgements}
	We thank Dr.~Kouta Sekine (Toyo University, Japan) for his helpful advice.
	We also express our gratitude to anonymous referees for insightful comments.
	This work was supported by CREST, JST Grant Number JPMJCR14D4;
	and by JSPS KAKENHI Grant Number JP19K14601.
	\bibliographystyle{elsarticle-num}
	
\end{document}